\documentclass[amscd,amssymb,verbatim,12pt]{amsart}
\usepackage{epsfig}
\usepackage{graphicx}

\newif\ifAMS
\IfFileExists{amssymb.sty}
  {\AMStrue\usepackage{amssymb}}
  {\usepackage{latexsym}}
  \usepackage{enumerate}
\theoremstyle{plain}

\newcommand{\kf}[1]{\marginpar{\tiny #1 --kf}}

\newtheorem*{condition}{Condition}

\newtheorem{Thm}{Theorem}[section]
\newtheorem{Cor}[Thm]{Corollary}
\newtheorem{Lem}[Thm]{Lemma}
\newtheorem{Prop}[Thm]{Proposition}

\theoremstyle{definition}
\newtheorem{Def}[Thm]{Definition}

\theoremstyle{remark}

\newtheorem{Rems}[Thm]{Remarks}

\DeclareMathOperator{\length}{length}
\DeclareMathOperator{\diam}{diam}

\newcommand{\interior}{^{ \kern-5pt ^\circ}}
\newcommand {\bd}{\partial}

\begin{document}
\title
{A coarse-geometry characterization of cacti}

\author
{Koji Fujiwara }
\email{kfujiwara@math.kyoto-u.ac.jp}
\address{Department of Mathematics, Kyoto University,
Kyoto, 606-8502, Japan}

\thanks{
The first author is
    supported in part by Grant-in-Aid for Scientific Research
    (No. 20H00114).}
\author
{Panos Papasoglu }

%\subjclass{53C23}

\email {} \email {papazoglou@maths.ox.ac.uk}

\address
{Mathematical Institute, University of Oxford, 24-29 St Giles',
Oxford, OX1 3LB, U.K.  }

\address
{ }

\begin{abstract} 
We give a quasi-isometric characterization of cacti, which 
is similar to Manning's characterization of quasi-trees by
the bottleneck property.
We also give another quasi-isometric characterization of cacti
using fat theta curves. 

\end{abstract}
\maketitle

\section{Introduction}

A \textit{cactus} is a connected graph such that any two cycles intersect at most one point.
Cacti play an important role in graph theory, in networks and computer science and in geometric group theory.
For example Dinitz-Karzanov-Lomonosov showed that one can encode the min-cuts of a graph by a cactus
\cite{DKL} (see also\cite{KT}, sections 7.4,7.5). This result has algorithmic applications to networks
where one would like to improve the edge connectivity of a graph (see \cite{NGM}). It turns out that one can also
encode the ends of infinite graphs by cacti (see \cite{EP}, \cite{Wi}) thereby deducing the classical Stallings' ends 
theorem in geometric group theory. 

Manning \cite{Ma} gave a coarse-geometry characterization of trees that has several interesting applications
in geometric group theory (see e.g. \cite{BBF}). 

One can think of cacti as the `simplest' graphs after trees. In this paper we study the coarse geometry of cacti. 
Our motivation comes from studying planar metrics (\cite{FP}). It turns out that the spheres of such metrics resemble cacti.
In this paper we make this precise by showing that such spheres are in fact quasi-isometric to cacti.
Our results apply more generally to geodesic metric spaces rather than graphs. We will call \textit{cactus} a geodesic metric space $C$ such that any two distinct simple closed
curves in $C$ intersect at at most one point. 
We note that our  notion of cactus generalizes the classical graph theoretic notion
in a similar way as $\mathbb R$-trees generalize trees.

Bonamy-Bousquet-Esperet-Groenland-Liu-Pirot-Scott, 
 \cite{BBEGLPS} showed that graphs with no $q$-fat $K_{2,p}$-minor (for some $p,q)$ have asymptotic dimension bounded by 1.
In this paper we show that if a graph has no $q$-fat $K_{2,3}$ minor then it is in fact quasi-isometric to a cactus
combining Corollary \ref{main} and Proposition \ref{prop-char-cactus}.
Further in \cite{BBEGLPS} the authors ask whether if a graph $\Gamma $ has no $q$-fat $H$-minor for some finite graph $H$ then
$asdim \Gamma <\infty $. Georgakopoulos-Papasoglu
\cite{GP}   made the following conjecture that would imply this: If a graph has no 
$q$-fat $H$-minor (for some $q$) then it is quasi-isometric to a graph with no $H$-minor. Our results confirm this conjecture for
$H=K_{2,3}$. This conjecture (if true) would show that the deep theory of minor free graphs generalizes to
$q$-fat minor free graphs. This would be interesting geometrically too as $q$-fat minor free graphs are invariant under quasi-isometries while
minor free graphs are not.

\subsection{Main results}

We proceed now to give our geometric characterization of
cacti up to quasi-isometry.

Let $X$ be a geodesic metric space. 
If $A$ is a subset $X$ and $m\in \mathbb R$ we set $$N_m(A)=\{x\in X:d(x,A)<m\}.$$

We consider the following condition:

\begin{condition}[$\sharp,m $]\label{condition.sharp}

  for any $x,y \in X$ with $d(x,y) \ge 10m$,
there exist $a,b$ with 
$d(x,a \cup b), d(y, a\cup b) \ge 4m$
and $x,y$ lie in distinct components of $X \backslash N_m(a\cup b)$.
\end{condition}
We may simply write  $(\sharp)$  without specifying the constant $m$.

Our main result is:

\begin{Cor}[Corollary \ref{bottlecor}, cf. Theorem \ref{bottle}]\label{main}
Let $X$ be a geodesic metric space.
$X$ is quasi-isometric to a cactus
if and only if it satisfies $(\sharp)$.

Moreover the quasi-isometry-constants depend only on $m$ in $(\sharp)$.
\end{Cor}

The hard direction is to show $(\sharp)$  implies
quasi-isometric to a cactus, which takes
most of the paper.

\subsection{Bottleneck property and quasi-trees}
The condition $(\sharp)$ is 
similar in spirit to the bottleneck property
that appears in 
Manning's Lemma \cite{Ma}: A geodesic metric space $X$ satisfies
{\em the bottleneck criterion} if there exists $\Delta \geq 0$ such that for any two points $x,y \in X$ the midpoint $z$ of a geodesic between $x$ and $y$ satisfies the property that any path from $x$ to $y$
intersects the $\Delta$-ball centered at $z$. Manning showed that this is equivalent to $X$ being a
quasi-tree. 
By definition, a geodesic space is called a {\em quasi-tree}
if it is quasi-isometric to a simplicial tree. 
His condition is rephrased as follows. We will not use this result in this paper.

\begin{Prop}[version of Manning lemma]
$X$ is a quasi-tree if and only if there exists $k$ such that for 
any $x,y \in X$ with $d(x,y) \ge M=100k$, 
there is $c\in X$ with $d(c,x), d(c,y) \ge 2k$ such that $x,y$ lie in distinct components
of $X \backslash B_c(k)$.
\end{Prop}
\proof
We apply Manning's lemma for $\Delta =100k$. Let $x,y\in X$ and let $z$ be the midpoint of the geodesic between $x,y$.
Let $x'\ne y'$ on this geodesic with $d(x',z)=d(y',z)=100k$. There is some point $c$ on the geodesic between $x',y'$ such that
$B_c(k)$ separates $x',y'$. But then $B_z(100k)$ contains $B_c(k)$ so it separates $x,y$.
\qed

Kerr showed that 
a geodesic metric space is a quasi-tree if and only if
it is $(1,C)$-quasi-isometric to a simplicial tree 
for some $C \ge 0$, \cite[Theorem 1.3]{Ke}.
It would be interesting to know if $X$ is a quasi-tree if 
and only if 
it is $(1,C)$-quasi-isometric to a cactus for some $C\ge0$.

\subsection{Background} 
This paper is a companion paper to \cite{FP}.
We explain the connection. 
The main result of \cite{FP} is:
\begin{Thm}
Let $P$ be a geodesic space that is homeomorphic to $\Bbb R^2$. Then the asymptotic dimension of $P$ is at most three, uniformly.

\end{Thm}

We review the strategy of the argument:
first we fix a point $e\in P$ and a constant 
$L$, then look at annuli, $A_n$,  that are sets of points whose distance 
from $e$ are between $nL$ and $(n+1)L$ for $n \in \Bbb N$.
We are lead to study $A_n$, and realized that those $A_n$ look like cacti. We did not 
show they are quasi-isometric to cacti, but 
we managed to show that the connected components of them are ``coarse cacti'' in the 
following sense \cite[Lemma 4.1]{FP}:

\begin{Def}[$M$-coarse cactus]
Let $X$ be a geodesic metric space. If there is an $M>0$ such that $X$ has no $M$-fat theta curves then we say that
$X$ is an $M$-\textit{coarse cactus} or simply a \textit{coarse cactus}.
\end{Def}

See Definition \ref{def-fat-theta} for the definition of fat theta curves. 

It is pretty easy  to show (see \cite{FP}):
\begin{Prop}\label{cactus.asdim1}
A cactus $C$ has $asdim \le 1$.
Moreover, $asdim \le 1$ uniformly over all cacti. 

\end{Prop}

The argument generalizes to show (\cite[Theorem 3.2]{FP}):
\begin{Thm}\label{quasi.cactus}
If $C$ is an $M$-coarse cactus then ${\rm asdim}\, C \le 1$.
Moreover, it is uniform with $M$ fixed. 
\end{Thm}

As a summary we \cite{FP} showed that those $A_n$ have asymptotic dimension 
at most $1$, uniformly.
Then, some general idea for dimension theory applies and 
it is easy for us to show Theorem 1.3.
In fact, since those $A_n$ are aligned along 
the radical direction from the base point, it is 
natural to expect the bound for $P$
is 2, but we left it as a question. 
By now, the question is solved positively and  as a result
the bound for $P$ is 2 (\cite{JL}, \cite{BBEGLPS}), which is optimal. 

We go back to the connection to the present paper. 
As the name suggests, 
we suspected that 
if $X$ is an ($M$-)coarse cactus then it is in fact quasi-isometric to a cactus,
and the quasi-isometry constants depend only on $M$.
Note that the converse implication is easy.
Combining Corollary \ref{main} and Proposition \ref{prop-char-cactus}
we will confirm the speculation (Lemma \ref{lem-annuli}).

\subsection{Fat theta curves and $(\sharp)$}

To discuss the relation between the no-fat-theta curve condition 
and $(\sharp)$, we recall some definitions precisely  from \cite{FP}.
Let $X$ be a geodesic metric space. Let $\Theta $ be a unit circle in the plane together with a diameter.
We denote by $x,y$ the endpoints of the diameter and by $q_1,q_2,q_3$ the 3 arcs joining them (ie the closures of the connected components of $\Theta \setminus \{x,y\}$).
A \textit{theta-curve} in $X$ is a continuous map $f:\Theta \to X$. Let $p_i=f(q_i),\, i=1,2,3,\, a=f(x),b=f(y)$.

We recall a definition from \cite{FP}.
\begin{Def}\label{def-fat-theta} [$M$-fat theta curve]
Suppose $M>0$. 
A theta curve is $M$-\textit{fat} if 
there are arcs $\alpha _i,\beta _i\subseteq p_i,\, i=1,2,3$ where $a\in \alpha _i,b\in \beta _i$
so that 
the following hold:
\begin{enumerate}
\item
If $p_i'=p_i\setminus \alpha _i\cup \beta _i$ then $p_i'\ne \emptyset $ and for
any $i\ne j$ and any $t\in p_i',s\in p_j'$ we have $d(t,s)\geq M$.
\item
$p_i'\cap \alpha _j=\emptyset,\, p_i'\cap \beta _j=\emptyset $ for all $i,j$ (note by definition $p_i'$ is an open arc,
ie does not contain its endpoints).
\item
For any $t\in \alpha _1\cup \alpha _2\cup \alpha _3, s\in  \beta_1\cup \beta _2\cup \beta _3$, we have $d(t,s)\geq 2M$.
\end{enumerate}
We say that $a,b$ are the \textit{vertices} of the theta curve.
We say that the theta curve is \textit{embedded} if the map $f$ is injective.
We will often abuse notation and identify the theta curve with its image giving simply the arcs of the theta curve.
So we will denote the theta curve defined above by $\Theta (p_1,p_2,p_3)$.

\end{Def}

If a geodesic space contain an $M$-fat theta curve, then 
it contains an embedded $M$-fat theta curve, which is a subset
of the first theta curve \cite[Lemma 3.1]{FP}.

As an easy application of Corollary \ref{main} we show:

\begin{Lem}[$\sharp$ implies no fat theta]\label{equivalence}
If a geodesic space $X$ satisfies $(\sharp, m)$ for some $m>0$, then 
there is $M>0$ such that $X$ has no $M$-fat theta curve. Moreover, the constant 
$M$ depends only on $m$.
\end{Lem}

It turns out the converse is true.
\begin{Lem}[No fat theta implies $\sharp$] \label{fatissharp}
Let $X$ be a geodesic metric space. Assume that
there is $M>0$ such that $X$ has no $M$-fat theta curve. Then $X$ satisfies $(\sharp, m)$ for some $m>0$.

 Moreover, the constant 
$m$ depends only on $M$.
\end{Lem}

Combining the two lemmas, we obtain the following 
proposition, so that the 
no-fat-theta condition gives another characterization 
for a geodesic space to be quasi-isometric to a cactus. 

\begin{Prop}\label{prop-char-cactus}
Let $X$ be a geodesic space.
Then $X$ satisfies $(\sharp, m)$ for some $m>0$
if and only if there is $M$ such that $X$ has no $M$-fat theta curve. 
\end{Prop}

See Section \ref{sec-char-proof} for the proof of the two lemmas and the proposition.

\section{Spaces Quasi-isometric to cacti}

We will prove our characterization of cacti for graphs (Theorem \ref{bottle}) but it will be easy to see that applies to all geodesic metric spaces (Corollary \ref{main}).
The purpose of this section is to prove those two results. 
For the rest of the paper by graph we mean a connected graph unless we specify otherwise.

We state our main result. 
\begin{Thm}[Manning lemma for cacti]\label{bottle}
Let $X$ be a graph.
$X$ is (uniformly) quasi-isometric to a cactus
if and only if there exists $m>10$ such that the  condition
$(\sharp,m) $
is satisfied.
It is uniform in the sense that the quasi-isometry-constants depends
only on $m$.
\end{Thm}

The proof of this theorem is similar to Manning's characterization lemma for quasi-trees but it is quite
more involved as we need to associate a cactus $C$ to the graph $X$, and in our construction of $C$ we have
to choose some (`big') simple closed curves of $X$ that will be preserved to the cactus $C$ while others will be collapsed
to intervals.

\subsection{Preliminary lemmas}
We prepare several lemmas before we start proving 
the theorem. 

\begin{Def}[Geodesic circle]
A simple closed curve $\alpha $ in a geodesic metric space $X$ is called a \textit{geodesic circle} if it is the
image of a circle (with its length metric) under an isometric embedding.
If $\alpha $ is a circle and $x,y\in \alpha $ we denote by $\overline{xy}$ the shortest of the two arcs joining
$x,y$ in $\alpha $ (if both arcs have equal length then we denote by $\overline{xy}$ any of them).
\end{Def}

\begin{Def}[Filling]
If $\alpha $ is a simple closed curve in a graph $X$ we say that a graph morphism $f:D\to X$ is
is a \textit{filling} of $\alpha $ if:

1) $D$ is a finite graph embedded in $\mathbb R^2$.

2) If $U$ is the unbounded component of $\mathbb R^2\setminus D$
then $\bd U$ is a simple closed curve such that $f(\bd U)=\alpha $.

3) If $\gamma $ is the boundary of any bounded connected component of $\mathbb R^2\setminus D$
then $f(\gamma )$ is a geodesic circle in $X$.

\end{Def}

We will often abuse notation and consider a filling $D$ of a curve in $X$ as a subset of $X$.
We will call the bounded connected components of $\mathbb R^2\setminus D$ \textit{regions of} $D$.

\begin{Lem}\label{filling}
Any finite length simple closed curve $\alpha $ in a graph $X$ has a filling.

\end{Lem}
\proof If $\alpha $ is a geodesic circle then there is nothing to prove.
Otherwise pick two points on $\alpha $ such that $d_{\alpha }(x,y)-d(x,y)$ attains its maximum (where
by $d_{\alpha }$ we denote the length of the shortest subpath of $\alpha $ joining $x,y$). Join then
$x,y$ by a path of length $d(x,y)$. We may see this subdivision as a map from a planar graph $D$ to $X$.
Now repeat this by considering any bounded region of $\mathbb R^2\setminus D$ that has a boundary that is not a 
geodesic circle. Note that the maximal boundary length drops by at least one after finitely many steps, so this
procedure terminates producing a filling (similar to a van-Kampen diagram).

\qed

It will be convenient in what follows to fix a constant $$M\gg m ,$$
for example $M=10^{100}m$ works for all the following lemmas.

Most of our lemmas concern inequalities that are (obvious) equalities in the case of cacti. In
our case we get inequalities instead, involving an `error term' expressed as a multiple of $m$. 
The precise expression in $m$ is not important- one may think of all multiples of $m$ appearing in
the sequel as `negligible' quantities. The proof is similar in most cases: we show that if the
inequality does not hold then $\sharp$ is violated. The way we show the latter is by producing
a `theta' curve, the branch points of which give us the points that violate $\sharp$.
Finally any time we manipulate inequalities we just state `obvious' inequalities rather than optimal
ones.

Quantities that are multiples of $M$ on the other hand are non-negligible, for example geodesic circles in $X$
of length greater than $10M$ are certainly represented in the corresponding cactus $C$.

\begin{Lem}\label{geodcirc} Let $X$ be a graph that satisfies $\sharp $ and let $\alpha $ be a geodesic circle in $X$.
If $e$ is a point in $X$, $d(e,\alpha )=d(e,p)=R$ with $p\in \alpha $ and if $x\in \alpha $ then
$$d(e,x)\geq R+\length (\overline{xp})-100m .$$
\end{Lem}
\proof
We argue by contradiction- so we assume the inequality of the lemma does not hold.
Consider a shortest path $\gamma =[x,e]$. Let $y $ be the last point on $\gamma $ (starting from $x$) such that $d(y,\alpha )=10m$.
Let $q\in \alpha $ such that $d(q,y)=10m$. We claim that $d(p,q)>50m$. If not then 
$$d(y,e)\geq R-60m$$ so
$$d(x,y)\leq \length (\overline{xp})-40m $$ and
$$d(q,x)\leq \length (\overline{xp})-30m $$ contradicting the fact that $\alpha $ is a geodesic circle.

\begin{figure}[htbp]
\hspace*{-3.3cm}     
\begin{center}
                                                      
   \includegraphics[scale=0.600]{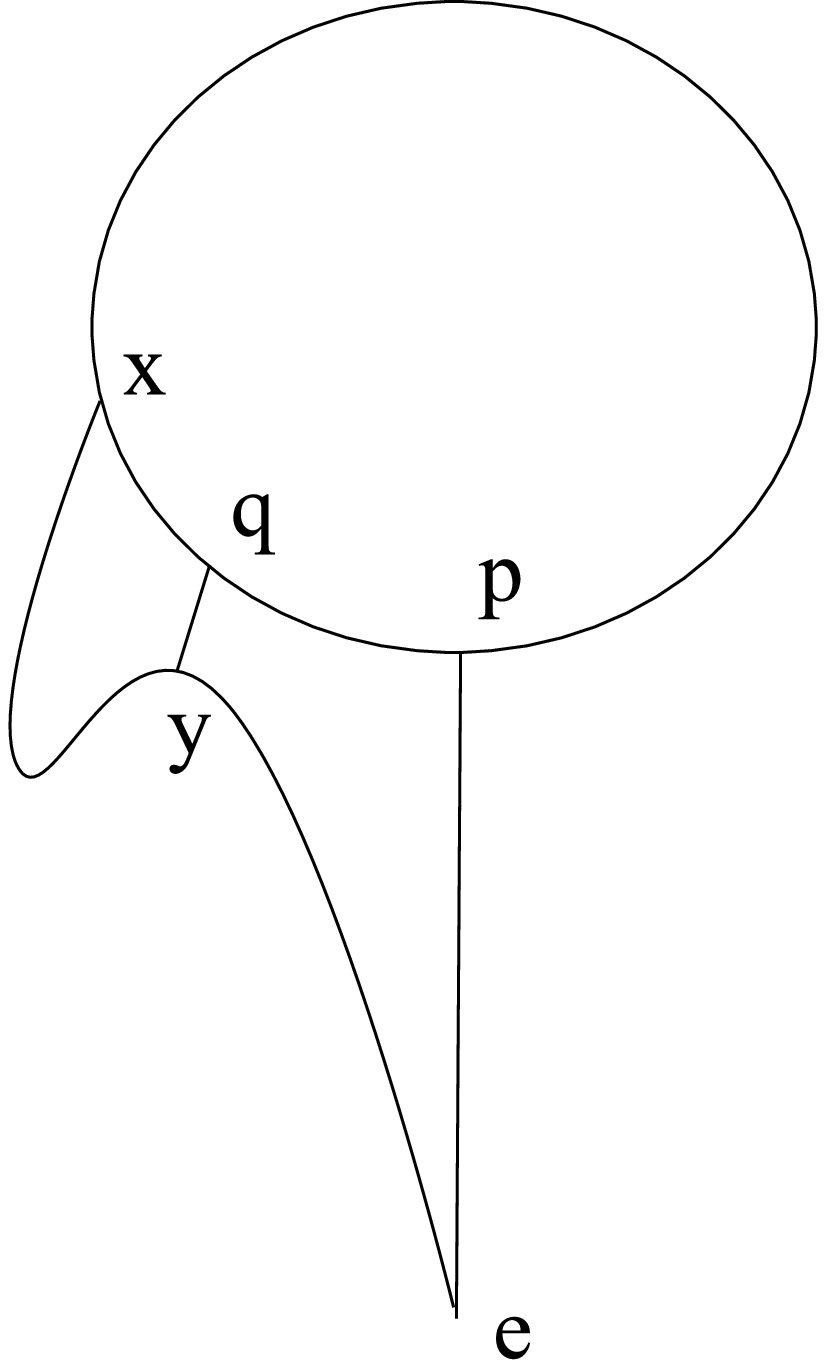}%
   \end{center}

\caption{}
  \label{}
\end{figure}

We claim now that condition $\sharp $ fails for $p,q$. Indeed since $p,q$ lie on $\alpha $ if $B_1=B(a_1,m),B_2=B(a_2,m)$ are two balls
of radius $m$ separating $p,q$ as in condition $\sharp $  then each one of $B_1$, $B_2$ should intersect one of the arcs on $\alpha $ joining $p,q$.
However we claim that neither $B_1$ nor $B_2$ intersects the arc $[q,y]\cup [y,e]\cup [e,p]$ joining $p,q$. Indeed if say $B_1$
intersects $[q,y]$ at $z$ then $d(z,q)\geq 3m-1$ by $\sharp $ and since $\alpha $ is geodesic but then 
$$d(y,\alpha )\leq d(y,z)+d(z,a_1)+m<10m $$ a contradiction. Similarly we arrive at the inequality $d(x,\alpha )<R$ if $B_1$ intersects $[x,p]$ and clearly
$B_1$ does not intersect $[e,y]$. Therefore $\sharp $ does not hold for $p,q$, a contradiction.

\qed

\begin{Lem}\label{circles-close} Let $X$ be a graph that satisfies $\sharp $ and let $S_1,S_2$ be two geodesic circles in $X$ with $$\length (S_1),\length (S_2)\geq M/10$$ and such that
for some $s\in S_1$, $$d(s,S_2)\geq M/100.$$
If $a,b\in S_1$ are such that $d(a,S_2), d(b,S_2)\leq 10m$ then $d(a,b)\leq 30m$.

\end{Lem}
\proof
We argue by contradiction.
Let $w_1,w_2:I\to S_1$ be the two parametrisations by arc length of $S_1$ starting from $s$.
So, say, $w_1$ goes first through $a$ while $w_2$ goes first through $b$. Let's say that $a=w_1(t_1), b=w_2(t_2)$.
Let $t$ be minimal such that $t\leq t_1$ and $d(w_1(t),S_2)=10m$.
Let $t'$ be minimal such that $t'\leq t_2$ and $d(w_2(t'),S_2)=10m$.
We set $c_1=w_1(t), d_1=w_2(t')$. Let $c_2,d_2$ be the closest points to $c_1,d_1$ on $S_2$ respectively.

We claim now that $c_2,d_2$ violate condition $\sharp $. Indeed since $d(a,b)> 30m$ we have that $d(c_2,d_2)\geq 10m$.
There are two distinct arcs, say $\alpha _1, \alpha _2$, on $S_2$ joining $c_2,d_2$. If $\alpha $ is the arc on $S_1$ joining $c_1,d_1$ which contains $s$
then we have the arc $$\alpha _3 =[c_2,c_1]\cup \alpha \cup [d_1,d_2]$$
joining $c_2,d_2$ as well. Let $B=B_z(m)$ be a ball with $d(z,\{c_2,d_2\})\geq 4m$. Since $S_2$ is 1-geodesic $B$ can not intersect both $\alpha _1, \alpha _2$.
By definition of $\alpha $, $B$ does not intersect $\alpha $ and some $\alpha _i, i=1,2$. If $B$ intersects some $\alpha _i$ (for $i=1,2$) and say $[c_1,c_2]$ then $d(c_1,S_2)\leq 9m$
contradicting the definition of $c_1$. Similarly $B$ does not intersect $\alpha _i$ for $i=1,2$ and $[c_1,c_2]$. So $B$ can not intersect two of the $\alpha _i$'s for $i=1,2,3$,
which is a contradiction.

\qed

\begin{Lem}\label{2circles} Let $X$ be a graph that satisfies $\sharp $ and let $S_1,S_2$ be two geodesic circles in $X$ with $$\length (S_1),\length (S_2)\geq M/10$$ and such that
for some $s\in S_1$, $$d(s,S_2)\geq M/100.$$

If $e,p$ are points in $S_1,S_2$ respectively such that
$$d(e,p)=d(S_1,S_2)=R$$ 
and if $x\in S_1,y\in S_2 $ then
$$d(x,y)\geq R+\length (\overline{xe})+\length (\overline{py})-1000m .$$
\end{Lem}
\proof
We argue by contradiction- so we assume the inequality of the lemma does not hold.
Consider a shortest path $\gamma =[x,y]$. If there is some point $z$ on $\gamma $ and some $z_1\in [e,p]$ such that $d(z,z_1)\leq 100m$
then $$d(x,z)\geq d(x,z_1)-100m $$ and by lemma \ref{geodcirc} $$d(x,z_1)\geq d(e,z_1)+\length (\overline{xe})-100m .$$
So $$d(x,z)\geq d(e,z_1)+\length (\overline{xe})-200m .$$
Similarly we get 
$$d(y,z)\geq d(p,z_1)+\length (\overline{py})-200m .$$
so $$d(x,y) \geq R+ \length (\overline{xe})+\length (\overline{py})-400m $$
and the lemma holds. So we may assume that $d(\gamma, [e,p])>100m$. We distinguish two cases:

\textit{Case 1.} $R>10m$. Let $y_1$ be the last point on $\gamma $ (starting from $x$) such that $d(y_1,S_1)=10m$ 
and let $x_1\in S_1$ such that $d(y_1,x_1)=10m$.
Then we claim that $e,x_1$ violate $\sharp $. 
Since $d(\gamma, [e,p])>100m$, $d(e,x_1)>10m$. Let $\alpha _1, \alpha _2$ be the two distinct arcs on $S_1$ joining $e,x_1$.

 We consider now the following arc joining $e,x_1$:

$$\alpha _3=[e,p]\cup \overline{py}\cup [y_1,y] \cup [x_1,y_1] .$$
Let $B=B_z(m)$ be a ball with $d(z,\{e,x_1\})\geq 4m$. Since $S_1$ is 1-geodesic $B$ can not intersect both $\alpha _1, \alpha _2$.
Since $R>10m$, $B$ does not intersect $\overline{py}$ and some $\alpha _i, i=1,2$. If $B$ intersects some $\alpha _i$ (for $i=1,2$) and $[e,p]$ then $d(S_1,S_2)\leq R-m$
contradicting our definition of $R$. Similarly $B$ does not intersect some $\alpha _i$ for $i=1,2$ and $[x_1,y_1]$ as we would have $d(y_1,S_1)\leq 9m$. Finally
$B$ does not intersect some $\alpha _i$ for $i=1,2$ and $[y_1,y]$ as this would contradict the definition of $y_1$.
So $B$ can not intersect two of the $\alpha _i$'s for $i=1,2,3$,
showing that $e,x_1$ violate $\sharp $.

\textit{Case 2.} $R\leq 10m$. Let $y_1$ be the first point on $\gamma $ (starting from $x$) such that $d(y_1,S_2)=10m$ 
and let $y_2\in S_2$ such that  $d(y_1,y_2)=10m$.
Consider the largest subarc $\overline {ab}$ of $S_1$ containing $e$ such that
$a,b$ are at distance $\leq 10m$ from $S_2$. By lemma \ref{circles-close} $d(a,b)\leq 30m$. 
Let $a_1,b_1\in S_2$ such that $$d(a,a_1)=d(b,b_1)=10m .$$
We claim that $a_1,y_2$ violate $\sharp $. 

Since $d(\gamma, [e,p])>100m$, and (by lemma \ref{circles-close}) $d(a_1,b_1)<30m$ clearly $d(a_1,y_2)>10m$.
Let $\alpha _1, \alpha _2$ be the two distinct arcs on $S_2$ joining $a_1,y_2$.
Let $\alpha $ be the arc of $S_1$ joining $a,x$ which does not contain $b$. We have the following arc joining $a_1,y_2$:
$$\alpha _3= [a_1,a]\cup \alpha \cup [x,y_1]\cup [y_1,y_2].$$
Let $B=B_z(m)$ be a ball with $d(z,\{a_1,y_2\})\geq 4m$.  Since $S_1$ is 1-geodesic $B$ can not intersect both $\alpha _1, \alpha _2$.
By the definition of $a,b$ and $y_1$ $B$ can not intersect some $\alpha _i$ and $\alpha $ or $[x,y_1]$ for $i=1,2$.
Finally if for some $i=1,2$ $B$ intersects $[a,a_1]$ or $[y_1,y_2]$ we have $d(a,S_2)\leq 9m$ (in the first case) or 
$d(y_1,S_2)\leq 9m$ (in the second case) contradicting the definition of $a$ or of $y_1$.
So $B$ can not intersect two of the $\alpha _i$'s for $i=1,2,3$,
showing that $a_1,y_2$ violate $\sharp $.

\qed

\begin{Lem}\label{circle} Let $X$ be a graph that satisfies $\sharp $ and let $Y\subseteq X$ be a connected subgraph of
$X$. Let $Z$ be a connected component of $X\setminus N_M(Y)$.
If  there are $a,b\in \bd Z$ with $d(a,b)\geq M/10$ then there is a geodesic circle $S$ in $X$ such that $a,b\in N_{100m}(S)$.
\end{Lem}
We will use this lemma to construct a cactus inductively, 
namely, apply the lemma to a ``node'' (a point or a circle) of a cactus as $Y$, then obtain a geodesic circle, which is a new node. 
\proof
Let $\alpha _1,\beta _1$ be shortest paths joining $a,b$ to $Y$ respectively.
Let $\gamma _1$ be a path in $X\setminus N_{M}(Y)$ joining $a,b$. 
If there is a point on $\alpha _1$ at distance $<100m $ from $\beta _1$ we set $a_1$ to be the first point (starting from $a$)
on $\alpha _1$ such that $d(a_1,\beta )=100m$ and we let $\gamma _2$ be a path of length $100m$ joining $a_1$ to $b_1\in \beta _1$.
We set then $\alpha $ to be the subarc of $\alpha _1$ with endpoints $a,a_1$ and $\beta $ to be the subarc of $\beta _1$
with endpoints $b,b_1$.

Otherwise we set $a_1,b_1$ to be respectively the endpoints of $\alpha _1,\beta _1$, we set $\alpha =\alpha _1, \beta = \beta _1$
and we take $\gamma _2$ to be a path in $Y$ joining $a_1,b_1$.

\begin{figure}[htbp]
\hspace*{-3.3cm}     
\begin{center}
                                                      
   \includegraphics[scale=0.700]{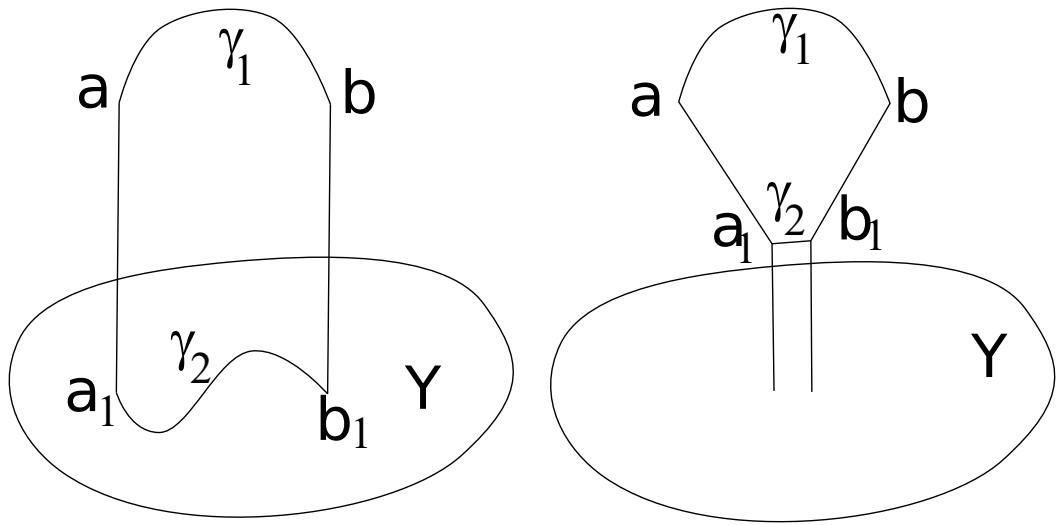}%
   \end{center}

\caption{Definition of $a_1,b_1,\gamma _2$}
  \label{}
\end{figure}

Let $$w=\alpha  \cup \gamma _1 \cup \beta \cup \gamma _2$$ and let $D$ be a filling of $w $. We abuse notation
and think $D$ as a subset of $X$. Let $a_2,b_2$ in $\alpha,\beta $ respectively such that $$d(a,a_2)=d(b,b_2)=4m$$
and $\gamma _3$ be a simple curve in $D\setminus N_{4m}(\alpha )$ joining points $a_2,b_2$ such that the component of 
$D\setminus \gamma _3$ containing $\gamma _1$ contains the least possible number of regions. We note that such a curve exists since there is a subarc
of $w$ joining $a_2,b_2$ in $D\setminus N_{4m}(\alpha )$ and if we have two such curves joining $a_2,b_2$ in
$D\setminus N_{4m}(\alpha )$ that cross then there is another curve joining $a_2,b_2$ which encloses less regions than either of the two.
We distinguish now some cases:\smallskip

\textit{Case 1}. $\gamma _3\cap N_{10m}(\gamma _2 )=\emptyset $. Let $x_1$ be the last point on $\gamma _3$ at distance $4m$ from $\alpha $ and
let $y_1$ be the first point on $\gamma _3$ at distance $4m$ from $\beta $. Let $x\in \alpha , y\in \beta $ such that 
$$d(x,x_1)=d(y,y_1)=4m .$$
We claim that condition $\sharp $ does not hold for $x,y$. Indeed consider the following three paths joining $x,y$: there are two paths joining them
on the simple closed curve $w$ and there is a path formed by the geodesic from $x$ to $x_1$ followed by the subpath of $\gamma _3$ with endpoints $x_1,y_1$
followed by the geodesic from $y_1$ to $y$. Using the triangle inequality one sees easily that for any $z\in X$ with $d(z,x),d(z,y)\geq 4m $
the ball $B_z(m)$ intersects at most one of these 3 arcs. Indeed if $B_z(m)$ intersects, say, the geodesic joining $x,x_1$ at $z_1$ and $\alpha $ then, as $d(z,x)\geq 4m$, by the triangle inequality $d(x_1,z_1)\leq m$ so
$d(x_1,\alpha )\leq 3m$ which contradicts our definition of $x_1$. The other cases are similar to this.

\begin{figure}[htbp]
\hspace*{-3.3cm}     
\begin{center}
                                                      
   \includegraphics[scale=0.600]{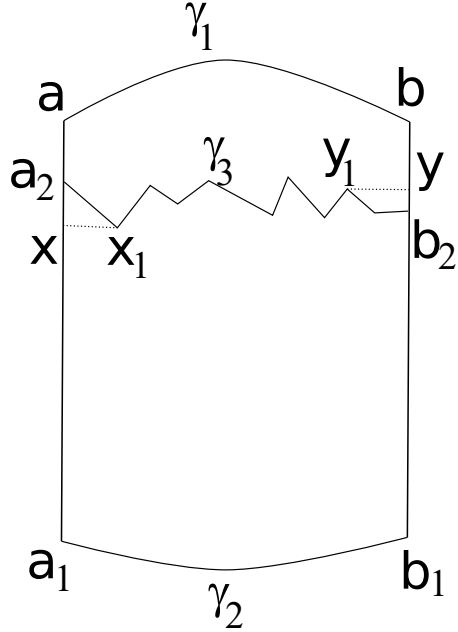}%
   \end{center}

\caption{Case 1: $\gamma _3$ is `far' from $\gamma _2$.}
  \label{}
\end{figure}

 \smallskip

\textit{Case 2}. $\gamma _3\cap N_{10m}(\gamma _2 )\ne \emptyset $. Let $x_1\in \gamma _3$ with $d(x_1,\gamma _2)\leq 10m$ and let $S$ be a geodesic circle in $D$
containing a non-trivial segment that contains $x_1$. We note that if $d(S,\gamma _1)>4m$ then we may replace $S\cap \gamma _3$ by the arc joining its endpoints
which does not intersect $\gamma _3$ in its interior contradicting our choice of $\gamma _3$ (as the new curve does not enclose the region bounded by $S$).
It follows that there is some $y_1\in \gamma _1$ such that $d(y_1,S)\leq 4m$. 

\begin{figure}[htbp]
\hspace*{-3.3cm}     
\begin{center}
                                                      
   \includegraphics[scale=0.600]{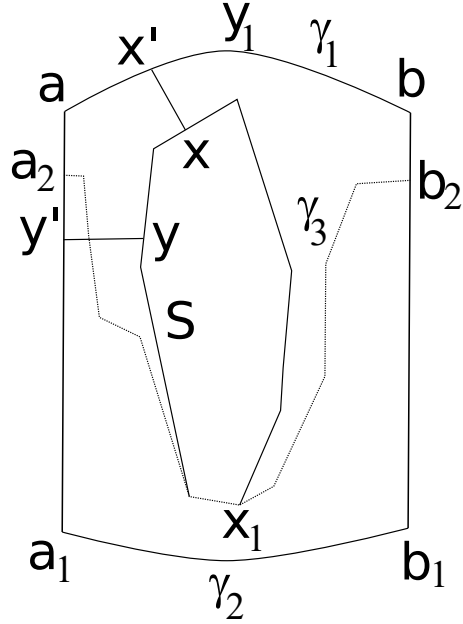}%
   \end{center}

\caption{Case 2: $\gamma _3$ is `close' to $\gamma _2$.}
  \label{}
\end{figure}

Clearly the lemma is proved if both $a,b$ are at distance $\leq 100m$ from $S$. So, without loss of generality, we assume that $d(a,S)>100m$.
There are two ways to traverse $w$ starting from $a$, so let $w_1:I\to w, w_2:I\to w$ be the corresponding parametrizations, where, say, for $t$ close to $0$, $$w_1(t)\in \gamma _1, w_2(t)\in \alpha .$$ Let $t_1$ be minimal such that $d(w_1(t_1),S)=10m$. Since $d(y_1,S)\leq 4m$ clearly $x'=w_1(t_1)\in \gamma _1$.
Let $t_2$ be minimal such that $d(w_2(t_2),S)=10m$. Since $d(x_1,S)\leq 10m$ clearly $y'=w_2(t_2)\in \alpha \cup \gamma _2$.
Let $x,y$ be points on $S$ such that $d(x,x')=d(y,y')=10m$. 
Note now that $$100m\leq d(y,a)\leq 10m+d(y',a)\Rightarrow d(y',a)\geq 90m.$$ It follows that
$d(y,Y)\geq M-80m$ so $d(y,x')\geq 80m$. Therefore $$d(x,y)\geq 70m .$$
We claim that $x,y$ violate condition $\sharp $.
There are 2 distinct arcs on $S$, say $S_1,S_2$ joining $x,y$. We consider also the following arc joining them
$$\eta=[x,x']\cup w_1([0,t_1])\cup  w_2([0,t_2])\cup [y,y'].$$
We claim that a ball $B_z(m)$ such that $$d(z,\{x,y\})\geq 4m$$ intersects at most one of these arcs. Since $S$ is a geodesic circle $B_z(m)$ can not intersect
both $S_1,S_2$. By the definition of $t_1,t_2$, $B_z(m)$ can not intersect one of the $S_i$'s and $ w_1([0,t_1])\cup  w_2([0,t_2])$. Let's say that
$B_z(m)$ intersects $S_1$ at $x_2$ and $[x,x']$ at $x_3$. Then, by the triangle inequality, $d(x,x_3)\geq 3m$. Since $d(x_2,x_3)\leq 2m$ it follows that
$d(x_2,x')\leq 9m$, which contradicts our hypothesis that $(x',S)=10m$. Clearly the same argument applies for $S_2$ or for $[y,y']$, so $\sharp $ is violated by $x,y$.
This finishes the proof of the lemma.

\qed

We show now that any point on the boundary of a connected component as before `close' to the circle we constructed at lemma \ref{circle}:
\begin{Lem}\label{close} Let $X$ be a graph that satisfies $\sharp $ and let $Y\subseteq X$ be a connected subgraph of
$X$. Let $Z$ be a connected component of $X\setminus N_M(Y)$.
Suppose there are $a,b\in \bd Z$ with $d(a,b)\geq M/10$.
Then if $c$ is any point in $\bd Z$
we have $$d(c,\{a,b\})\leq \dfrac {M}{10}.$$
\end{Lem}
\proof
We argue by contradiction, so we assume that $d(a,c)\geq M/10$.  By lemma \ref{circle} there is a geodesic circle $S_1$ in $X$ such that $a,b\in N_{100m}(S_1)$, 
also, 
 there is a geodesic circle $S_2$ containing $a,c$ in its $100m$ neighborhood.
We consider shortest paths $[e,p],[e,q]$ from the base point to the geodesic circles $S_1,S_2$ respectively. There are points $a_1\in S_1$, $a_2\in S_2$ such that
$$d(a,a_1)\leq 100m, d(a,a_2)\leq 100m $$
so $d(a_1,a_2)\leq 200m .$ If $d(b,c)\leq M/10$ the lemma is proven, so we may assume $d(b,c)>M/10$. Let $b'\in S_1$ with $d(b,b')\leq 100m $. We distinguish now two cases. \smallskip

\textit{Case 1.} $d(p,S_2)\leq 200m$.
There are two ways to traverse $S_1$ starting from $b'$, so let $w_1:I\to S_1, w_2:I\to S_1$ be the corresponding parametrizations. Let $t_1,t_2$ be minimal such that 
$$d(w_1(t_1),S_2)=200m, d(w_2(t_2),S_2)=200m .$$

\begin{figure}[htbp]
\hspace*{-3.3cm}     
\begin{center}
                                                      
   \includegraphics[scale=0.600]{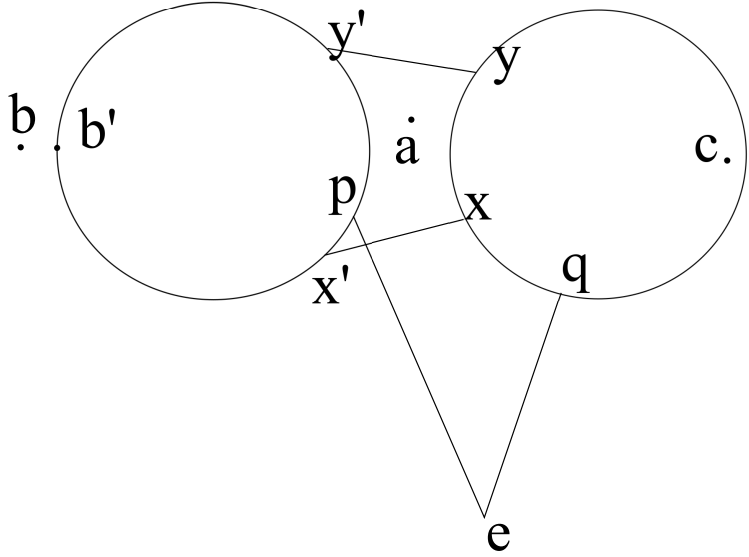}%
   \end{center}

\caption{Case 1: $p$ is `close' to $S_2$.}
  \label{}
\end{figure}

Let $x'=w_1(t_1), y'=w_2(t_2)$ and $x,y\in S_2$ such that
$$d(x',x)=d(y',y)=200m. $$
By lemma \ref{geodcirc} $d(p,a_1)\geq M/30$ so $d(x',y')>M/30$ and by the triangle inequality $d(x,y)>M/40>100m$.
We claim that $x,y$ violate condition $\sharp $.
Indeed we have the following 3 arcs joining $x,y$: there are two disjoint arcs on $S_2$ and if $\beta '$ is the arc of $S_1$ joining $x',y'$ which contains $b'$ we have also the arc
$$\beta =[x,x']\cup \beta '\cup [y,y'].$$
Let $B=B_z(m)$ be a ball with $d(z,\{x,y\})\geq 4m$. Then clearly $B$ can not intersect both $S_1,S_2$ and can not intersect any of the $S_i$ and $\beta '$. by the definition of $\beta '$.
Assume now that say $B$ intersects $[x,x']$ at $z_1$ and $S_1$ at $z_2$. Then $$d(z_1,x)\geq 3m, d(z_2,x)\geq 3m$$ and $d(z_1,z_2)\leq 2m $ so $d(z_2, x')<d(x,x')$
contradicting the definition of $x$. The other cases are similar so $x,y$ violate condition $\sharp $ in this case. \smallskip

\textit{Case 2.} $d(p,S_2)> 200m$. We distinguish two further cases:

\textit{Case 2a.} There is some $x'\in [e,p]$ such that $d(x',S_2)\leq 200m $. Without loss of generality
we assume that $x'$ is the last point on $[e,p]$ satisfying this property, so in particular $d(x',S_2)= 200m $. Let $x\in S_2$ such that
$d(x,x')=200m$. We pick now a parametrization $w:I\to S_1$ where $w(0)=b'$ and $a_1$ is reached before $p$. Let $t_1$ be minimal such that 
$$d(w(t_1),S_2)=200m,$$ we set as before $y'=w(t_1)$, $y\in S_2$ such that $d(y,y')=200m$. It follows by lemma \ref{geodcirc} that $d(x,y)\geq 100m$.

\begin{figure}[htbp]
\hspace*{-3.3cm}     
\begin{center}
                                                      
   \includegraphics[scale=0.600]{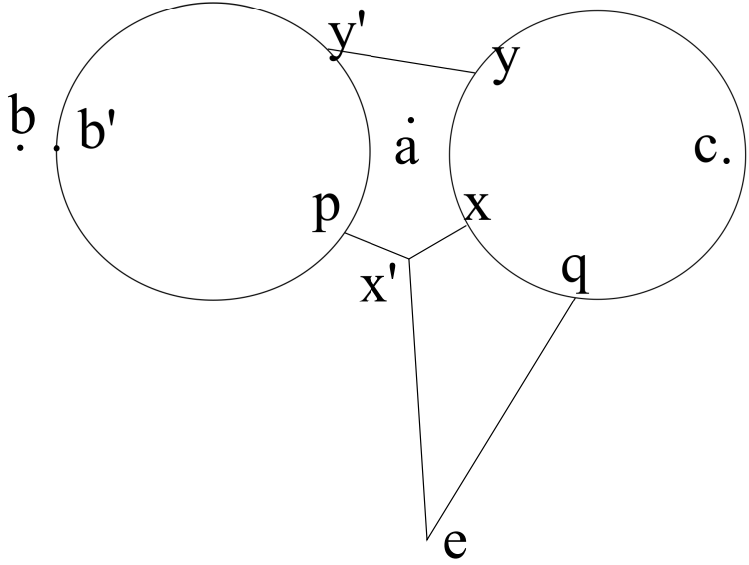}%
   \end{center}

\caption{Case 2a: Some point on $[e,p]$ is `close' to $S_2$.}
  \label{}
\end{figure}

We claim that $x,y$ violate condition $\sharp $.
Indeed we have the following 3 arcs joining $x,y$: there are two disjoint arcs on $S_2$ and if $\beta '$ is the arc on $S_1$ joining $p,y'$ passing through $b'$ we have also the arc
$$\beta =[x,x']\cup [x',p]\cup \beta '\cup [y,y'].$$
As is case 1 we see that a ball of radius $m$ can not intersect 2 of these arcs leading to a contradiction in this case.

\textit{Case 2b.} $d([e,p],S_2)>200m$. We define $y,y'$ as in case 2a. We claim that $q,y$ violate condition $\sharp $. Indeed by lemma \ref{geodcirc} that $d(q,y)\geq 100m$.
We have the following 3 arcs joining $p,y$: there are two disjoint arcs on $S_2$ and if $\beta '$ is the arc on $S_1$ joining $p,y'$ passing through $b'$ we have also the arc
$$\beta =[e,q]\cup [e,p]\cup \beta '\cup [y,y'].$$
As is case 1 we see that a ball of radius $m$ can not intersect 2 of these arcs leading to a contradiction in this final case too.

\qed
\subsection{Proof of Theorem \ref{bottle}}

We start the proof of Theorem \ref{bottle}.
\proof
%We return now to the proof of the proposition. 
We construct first inductively a cactus $C$ which we will show is quasi-isometric to $X$.
Let $e$ be a vertex of $X$ (that we think as basepoint). At step $k$ we will define a cactus $C_k$ and a 1-Lipschitz map $h_k:C_k\to X$.

To facilitate our construction we define also at each stage some subgraphs of $X$, so at step $k$ we define the subgraphs $Y^k_i,\, i\in I_k$.
Each $Y^k_i$ is associated to a node of $C_k$. We say that the $Y^k_i$ are the graphs of level $k$ of the construction. \smallskip

\textit{Step $0$}.
$C_0$ is the point $e$, which we call the node of level $0$.
Also, $e$ is the graph of level $0$ (in $X$).
Define $h_0(e)=e$.

\textit{Step $1$}.
We consider the connected components of $X\setminus N_M(e)$. Let $Y$ be such a component
and let $$\bd Y= Y \cap \overline {N_M(e)}.$$
We have 2 cases:

\textit{Case 1}. $\diam (\bd Y)<M/10$. Then we pick some $y\in \bd Y$. 

\textit{Case 2}. There are $a,b\in \bd Y$ with $\diam (\bd Y)=d(a,b)\geq M/10$ and applying lemma \ref{circle} there is a geodesic circle $S$
containing $a,b$ in its $100m$ neighborhood. In this case we have the following:

\begin{Lem} \label{contained} $\diam S\cap Y\geq M/20$. If  $Y'$ is any other component of $X\setminus N_M(e)$ then $\diam S\cap Y'\leq 100m$.
\end{Lem}

\proof
Let $a_1,b_1\in S$ be points in $S$ such that $$d(a,a_1)=d(a,S),\ \  \ d(b,b_1)=d(b,S).$$
Let $p\in S$ such that $d(e,p)=d(e,S)=R$. By lemma \ref{geodcirc} $$M+100m\geq d(e,a_1)\geq R+\length (\overline {pa_1})-100m.$$
and a similar inequality holds for $b_1$.
Since $$\length (\overline {a_1b_1})\geq \dfrac {M}{10}-200m$$
there is some point $x\in \overline {a_1b_1}$ such that $$d(x,e)\geq M+\dfrac {M}{30}$$
which shows that there is a subarc $\overline {a_2b_2}$ of $S$ of diameter $>M/20$ contained in a connected component $Y_1$ of $X\setminus N_M(e)$.
We may assume that $\overline {a_2b_2}$ is maximal with this property.

\begin{figure}[htbp]
\hspace*{-3.3cm}     
\begin{center}
                                                      
\includegraphics[scale=0.500]{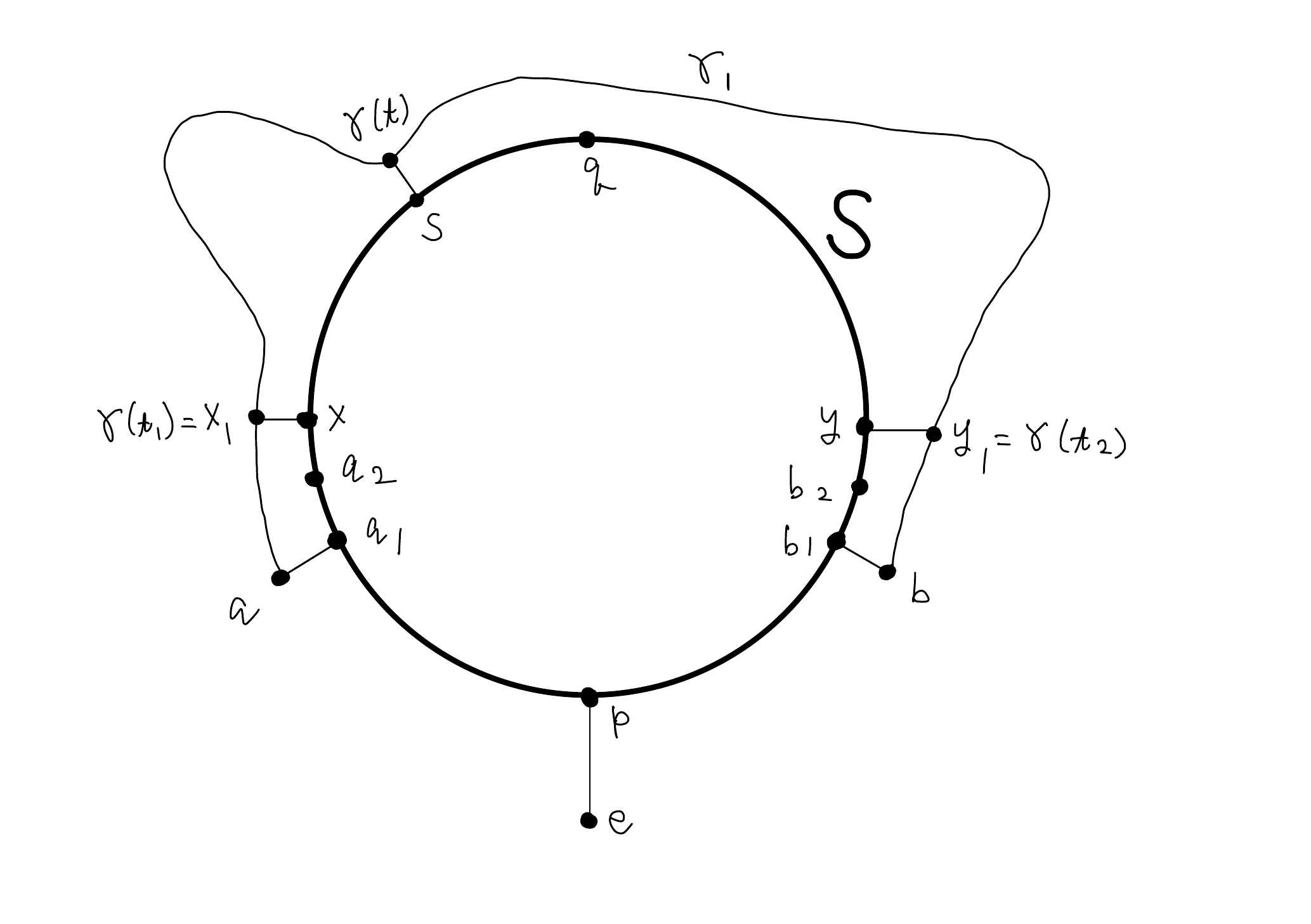}
  \end{center}
\caption{$d(a,a_1) \le 100m$, $d(x_1,a_1)=600m$, $d(x_1,x) \le 600m$, 
so $d(x,a) \le 1300m$. Similar for $b,b_1,y,y_1$.}
\label{contained.fig}
\end{figure}

We claim that $Y_1=Y$, which implies that $\diam S \cap Y \ge M/20$.
Assume that $Y_1\ne Y$. Let $\gamma _1$ be a path joining $a,b$ in $Y$.
See Figure \ref{contained.fig}. 
 Let
$$\gamma =[a_1,a]\cup \gamma _1\cup [b,b_1].$$
We parametrize $\gamma $ from $a_1$ to $b_1$ (so $\gamma (0)=a_1$). 
We note that if $d(\gamma (t),S)\leq 100m $ then either $d(\gamma (t), a_1)\leq 500m$
or $d(\gamma (t), b_1)\leq 500m$.
To see it, suppose not, and  let $s\in S$ be with $d(\gamma(t),s)=d(\gamma(t),S)$.
Then $d(s,a_1) >400m$ and $d(s,b_1) > 400m$.
By Lemma \ref{geodcirc}, either, $d(s,e) \le M - 200m$ (ie, $s$ gets
closer to $e$ compared to $a_1,b_1$), or 
$d(s,e) \ge M+ 200m$ (ie, $s$ gets farther from $e$ compared to $a_1,b_1$).
In the first case, we have $d(\gamma(t), e) \le M-100m$, so that 
$\gamma(t) \in N_M(e)$, impossible. In the second case, 
we have $s \in Y_1$, so that $Y$ and $Y_1$ are connected
by $[\gamma(t),s]$ outside of $N_M(e)$, impossible since $Y\not=Y_1$.

Let $t_1$ be maximal and $t_2$ minimal such that
$$d(\gamma (t_1),a_1) = 600m, d(\gamma (t_2),b_1)=600m.$$
Since $d(a_1,b_1) \ge M/20$, such $t_1,t_2$ exist and $t_1<t_2$.
Then we have $d(\{a_1,b_1\}, \gamma ([t_1,t_2])) \ge 600m$, so that 
$$d(S, \gamma ([t_1,t_2])) \ge 100m$$
as we said.

Let $$x_1=\gamma (t_1),\ \  y_1=\gamma (t_2)$$ and let $x,y\in S$ with
$$d(x,x_1)=d(S,x_1), d(y,y_1)=d(S,y_1), $$
which are at most $600m$. 
Then, by the triangle inequality, $d(x,y) \ge d(a,b)- 2600m \ge M/20$.

We claim that $x,y$ violate $\sharp $. Indeed, to see this, it suffices to consider the two arcs joining them in $S$
and the arc $$[x,x_1]\cup \gamma ([t_1,t_2])\cup [y_1,y] $$
since $d(S, \gamma ([t_1,t_2])) \ge 100m$.
It follows that $Y_1=Y$.

Lastly, 
let's say that $S$ intersects also another connected component $Y'$ of $X\setminus N_M(e)$.
Let $\overline {a_3b_3}$ be a maximal subarc of $S$ that intersects $Y'$.
$q$ necessarily lies in $\overline {a_2b_2}$. Hence $\overline {a_3b_3}$ is contained in either $\overline {pa_2}$ or in $\overline {pb_2}$
since $d(e,p) \le M-M/30$.
Let's say that $\overline {a_3b_3}$ is contained in $\overline {pa_2}$ and $b_3$ is closer to $a_2$. 
If $\length (\overline {a_3b_3})> 100m $ since

$$M=d(e,a_3)\leq R+\length (\overline {pa_3})$$ we have that

$$M=d(e,a_2)\geq R+\length (\overline {pa_2})-100m\geq R+\length (\overline {pa_3})+\length (\overline {a_3b_3})-100m>M$$
which is a contradiction.
\qed

\begin{figure}[htbp]
\hspace*{-3.3cm}     
\begin{center}
                                                      
\includegraphics[scale=0.400]{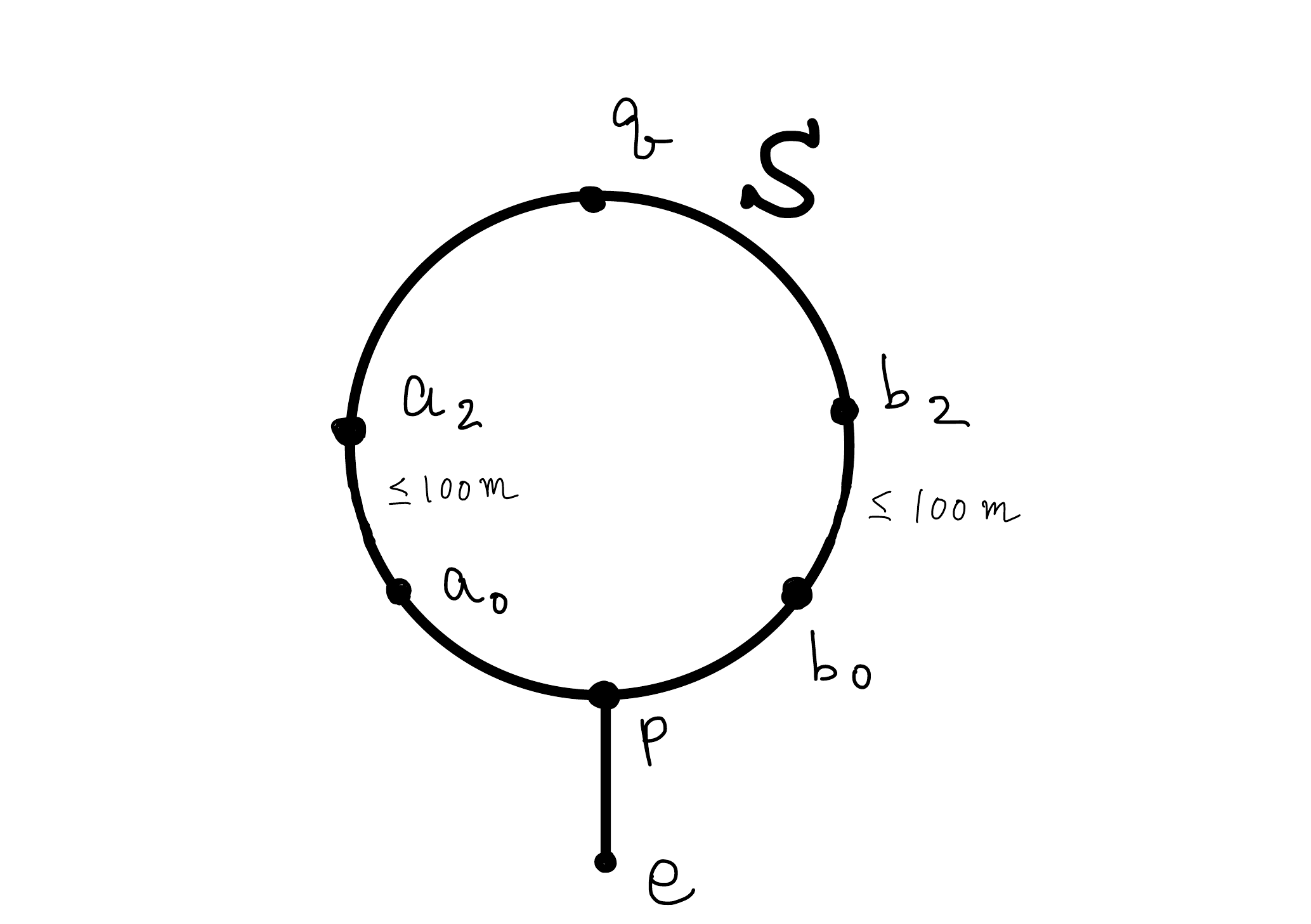}
  \end{center}
  
\caption{$\overline{a_2q} \cup \overline{qb_2}$ is a (longest) component 
of $S\cap Y$, and $\overline{a_0p}\cup \overline{pb_0}$
is a (longest) component of $S \cap N_M(e)$. Then $d(a_2,a_0), d(b_2,b_0)\le 100m.$}
\label{contained.remark.fig}
\end{figure}

\begin{Rems}\label{contained.remark}
%\kf{added}
Let $q \in S$ be the antipodal point from $p$. 
$S$ is diveided into two arcs, $A, B$, by $p$ and $q$.
Then
$d(e,q) \ge M+M/30$. This is by Lemma \ref{geodcirc}.
Now, let $Y'$ be the component of $X\backslash N_M(e)$ that 
contains $q$, and let $a_2 \in A, b_2 \in B$ such that 
$\overline{a_2q} \cup \overline{q b_2}$ is the component of $Y' \cap S$
that contains $q$. Note that $d(e,a_2)=d(e,b_2)=M$.
Note that by Lemma \ref{contained}, $Y=Y'$.

Also, note that $d(e,p)\le M-M/30$, so  let 
$a_0 \in A, b_0 \in B$ such that $\overline{a_0p} \cup \overline{pb_0}$ is the 
component of $S \cap N_M(e)$ that contains $p$. Note that $d(e,a_0)=d(e,b_0)=M$,
so that the length of $\overline{a_0p} \cup \overline{pb_0}$ is $\le 2M +200m$ by Lemma \ref{geodcirc}.
See Figure \ref{contained.remark.fig}. 

Now notice that $d(a_0,a_2), d(b_0,b_2) \le 100m$
by Lemma \ref{geodcirc}.

\end{Rems}

We go back to the construction. 
In the first case we join $e$ to $y$ by a segment of length $M$ while in the second case we join $e$ to a closest 
point of $e$ to $S$ by a segment of length equal to this distance. Note that in fact $e$ could lie on $S$, in this case the segment has 0 length (so we add nothing). We do this for all 
connected components and we obtain a cactus $C_1$. There is an obvious 1-Lipschitz map $h_1:C_1\to X$, we map $S,y$ as above by the identity map and the joining arcs to corresponding
geodesic segments. We say that $y,S$ are \textit{nodes} of level $1$ of the cactus. 

We index the connected components of $X\setminus N_M(e)$ by a set $I_1$, so at step $1$ we define the graphs $Y^1_i$ where $Y^1_i$ is a connected component of $X\setminus N_M(e)$.
Using this notation, by the previous discussion, if $$\diam (\bd Y^1_i)<M/10$$ then there is a node of the cactus $C_0$ that is a point which we denote by $y^1_i$ that lies in $\bd Y^1_i$,
while if $$\diam (\bd Y^1_i)\geq M/10$$ then by lemma \ref{contained} there is a unique node of $C_1$ which is a geodesic circle and which we denote by $S^1_i$ such that 
$$\diam (S^1_i\cap Y^1_i)\geq M/20 .$$ We note that there are two metrics on the subgraphs $Y^1_i$: the path metric of the subgraph and the induced metric from $X$. Even though
it is not very important for the argument we need to fix a metric to use, and we will be always using the \textit{induced metric} for subgraphs of $X$, unless we specify otherwise.

 \smallskip

\textit{Step $k+1$}. We consider the connected graphs $Y^k_i$ 
of level $k$ defined at step $k$. There are two cases in the level $k$, namely
the node is a point or a circle, then each case will have two cases
in the level $k+1$, ie, the node is a point or a circle. 

First, if
$$\diam (\bd Y^k_i)<M/10$$ then there
is a node of level $k$,
 $z=y^k_i$ of $C_k$ that is a point on $\bd Y^k_i$.

We consider the connected components of $Y^k_i\setminus N_M(z)$. Let $Y$ be such a component
and let $$\bd Y= Y \cap \overline{N_M(z)}.$$
We define $Y$ to be a graph of level $k+1$ of our construction.

We note that in fact $Y$ is a connected component of $X\setminus \overline{ N_M(z)}$ as well since $\diam (\bd Y^k_i)<M/10$, so that $\bd Y^k_i$ is contained
in $N_M(z)$.
(This is because the  only possibility that $Y$ becomes larger
as a connected component in $X\setminus \overline{ N_M(z)}$  is that $Y$ is connected to something
else at $\bd Y^k_i$. But this does
not happen since $\bd Y^k_i$ is contained in $N_M(z)$.)

If $\diam (\bd Y)<M/10$ then we pick some $y\in \bd Y$. Otherwise there are $a,b\in \bd Y$ with $\diam (\bd Y)=d(a,b)\geq M/10$ and applying lemma \ref{circle} to the connected component $Y$ (by considering the point $z$ as a subgraph in $X$, then 
$Y$ is a connected component of $X \backslash \overline{N_M(z)}$ as we noted)
there is a geodesic circle $S$
containing $a,b$ in its $100m$ neighborhood. 

%\kf{changed since the node is a point}
In the first case, we join the node $y$  to the node $z$ of $C_k$ by 
adding an edge of length $M=d(z,y)$ between $y$ and $z$.
In the case of $S$ we pick $s\in S$ such that $d(S,z)=d(s,z)$ and we add an edge of length $d(s,z)$ joining $s,z$. 

We note that lemma \ref{contained} applies  to $Y$ and $S$ in this case too since $Y$ is a connected component of $X\setminus \overline{ N_M(z)}$ as we noted. 
%\kf{added}

Second, if $$\diam (\bd Y^k_i)\geq M/10$$ then there
is a node of level $k$, $S'=S^k_i$ of $C_k$ that is a geodesic circle such that $\diam (S'\cap Y^k_i)\geq M/20$.

We consider the connected components of $Y^k_i\setminus N_M(S')$. Let $Y$ be such a component
and let $$\bd Y= Y \cap \overline{ N_M(S')}.$$

We note that, as before,  $Y$ is a connected component of $X\setminus \overline{ N_M(S')}$ as well since $\partial Y_i^k$ is contained in $N_{100m+M/10}(S')$
by Lemma \ref{close}.

%\kf{added.}

Similarly in this case we define $Y$ to be a graph of level $k+1$ of our construction.
As before, there are two cases. 
If $\diam (\bd Y)<M/10$ then we pick some $y\in \bd Y$. Otherwise there are $a,b\in \bd Y$ with $\diam (\bd Y)=d(a,b)\geq M/10$, then as in the previous case, applying lemma \ref{circle} there is a geodesic circle $S$
containing $a,b$ in its $100m$ neighborhood.  
Let $a_1,b_1\in S$ be points on $S$ such that $$d(a,a_1)=d(a,S),\ \  \ d(b,b_1)=d(b,S).$$ We have
$d(a,a_1),d(b,b_a) \le 100m$.

\begin{figure}[htbp]
\hspace*{-3.3cm}     
\begin{center}
                                                      
\includegraphics[scale=0.500]{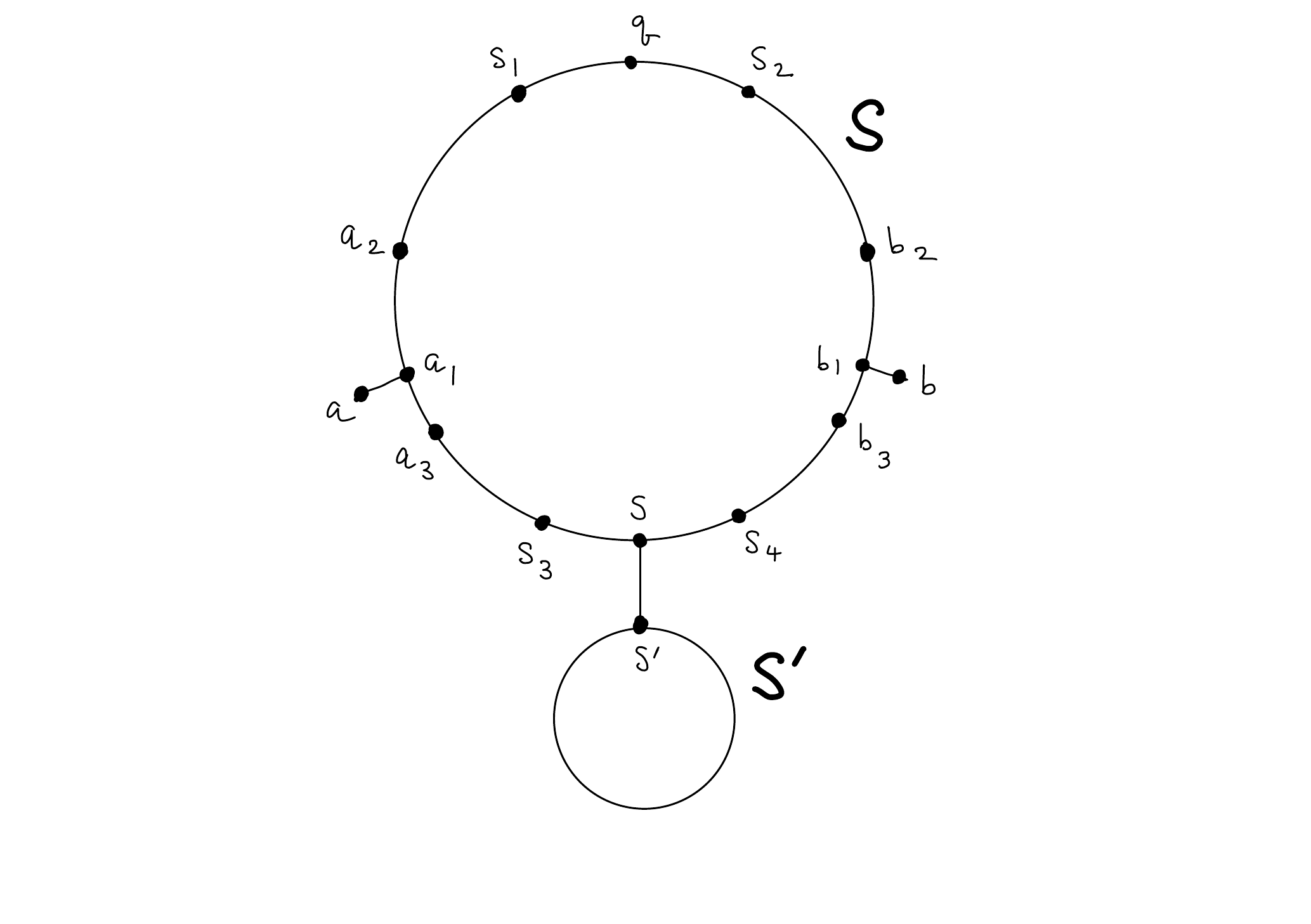}
  \end{center}
\caption{The left half of $S$ is $S_1$ and the right 
half is $S_2$. 
$\overline{a_2qb_2} \subset S \backslash N_M(S')$ and $\overline{a_3sb_3} \subset N_M(S').$}
\label{contained2.fig}
\end{figure}

\begin{Lem}\label{bigpart}
Suppose the nodes of level $k,k+1$ are circles $S,S'$.
Let $s\in S, s'\in S'$ be points with $d(S,S')=d(s,s')$.
Let $q\in S$ be the antipodal point from $s$.
Then
$$ d(q,S') >M, \, d(s,S')<M. $$

\end{Lem}
\proof
Let $S_1, S_2$ be the two arcs of $S$ with the endpoints $s,q$.
To argue by contradiction, 
suppose the first inequality fails, ie, $d(q,S') \le M$. 
Let $s_1\in S_1, s_2\in S_2$ be points with $s(q,s_i)=2000m$.
Since $d(s,q) \ge M/10$ such points exist.
Then by Lemma \ref{2circles}, $\overline{s_1s} \cup \overline{s_2s} \subset
N_{M-200m}(S')$.
This implies that $a_1,b_1 \in \overline{s_1s_2}$ since $d(a_1,S') \ge M-100m$ and $d(b_1,S') \ge M-100m$, so that 
$d(a_1,b_1) \le 4000m$, which contradicts $d(a_1,b_1) \ge M/10 - 200m > 4000m$.
See Figure \ref{contained2.fig}.

Suppose the second inequality fails, ie, $d(s,S') \ge M$.
Let $s_3 \in S_1, s_4 \in S_2$ with $d(s,s_i)=2000m$.
Then $\overline{s_3 q}\cup \overline{s_4 q} \subset S\backslash N_{M+200m}(S')$
as before by Lemma \ref{2circles}. So, $a_1,b_1 \in \overline{s_3s} \cup \overline{s_4 s}$, hence
$d(a_1,b_1) \le 4000m$, impossible. 
\qed

From now on, for example, instead of writing $\overline{s_1s} \cup
\overline{ss_2}$ we may write $\overline{s_1ss_2}$, which is maybe longer
than half of the circle. 

As a consequence of this lemma,
similarly as in lemma \ref{contained},  a `big part' of $S$ is contained in $Y$:
\begin{Lem}\label{contained2-2}
 $\diam S\cap Y\geq M/20$. If $Y'$ is any other component of $X\setminus N_M(S')$ then $\diam S\cap Y'\leq 2000m$.
\end{Lem}

We can prove this lemma in the same way as Lemma \ref{contained},
but we give a slightly different argument by proving 
something similar to Remark \ref{contained.remark} first. 
\proof
We use the notation of Lemma \ref{bigpart} and its proof. 
Also see Figure \ref{contained.fig}.

By Lemma \ref{bigpart}, $d(q,S')>M$, so let $a_2 \in S_1$ be such that 
$\overline{a_2q} \subset S\backslash N_M(S')$ and that $\overline{a_2 q}$
is maximal among such arcs. Similarly, let $b_2 \in S_1$
be such that $\overline{b_2q}$ is the maximal arc in $S\backslash N_M(S')$.
Note that $d(a_2, S')=d(b_2,S')=M$.

Also, by Lemma \ref{bigpart}, $d(s,S')<M$, let $a_3\in S_1$ be such that
$\overline{a_3s}$ is the maximal arc that is contained in the closure of 
$N_M(S')$ and $b_3 \in S_2$ such that 
$\overline{b_3s}$ is the maximal arc that is contained in the 
closure of $N_M(S')$.
Note that $d(a_3,S')=d(b_3,S')=M$.

Now, by Lemma \ref{2circles}, 
$d(a_2,a_3) \le 2000m, d(b_2,b_3) \le 2000m$.
Also, $d(a_1,a_2) \le 2100m, d(b_1,b_2) \le 2100m$.
It implies that 
$$d(a_2,b_2) \ge d(a_2,a_1) + d(a_1,b_1) +d(b_1,b_2) \ge
 M/10 -200m -4200m \ge M/20.$$
 
 Now let  $Y_1$ be the component of $X \backslash N_M(S')$
that contains $q$. Then it contains $\overline{a_2 q b_2}$.
Since $d(a_2,b_2) \ge M/20$, we have $\diam S \cap Y \ge M/20$.

We want to show $Y=Y_1$. The argument
is very similar to the one we showed $Y=Y_1$ in the 
proof of  Lemma \ref{contained}, so we omit it.

Lastly, if there is any other component $Y'$ than $Y$ that intersects $S$, then
$S\cap Y'$ must
be contained in $\overline{a_2a_3}$ or $\overline{b_2b_3}$, so that 
$\diam S \cap Y' \le 2000m$.
\qed

 We join $y$ (or $S$) to a node $S'$ of $C_k$ which is a shortest distance from $y$ (or $S$). So we pick $s'\in S'$ such that
$d(y,S)=d(y,s')$ and we add an edge of length $d(y,s')$ joining $y,s$. In the case of $S$ we pick $s\in S, s'\in S'$ such that $d(S,S')=d(s,s')$ and we add an edge of length $d(s,s')$ joining $s,s'$. We note that possibly $d(y,s')$ (or $d(S,S')$) is equal to 0.

We do this for all 
connected components and we obtain a cactus $C_{k+1}$. There is an obvious 1-Lipschitz map $h_{k+1}:C_{k+1}\to X$, we map $S,y$ as above by the identity map and the joining arcs to corresponding
geodesic segments. We note that any two points in $h_{k+1}(C_{k+1})$ can be joined by an arc of finite length in $h_{k+1}(C_{k+1})$. 

Finally we index all the graphs of level $k+1$ that we defined earlier by an index set $I_{k+1}$.

\smallskip

We set $$C=\bigcup _{k\in \mathbb Z}C_k .$$
%\kf{does it start with k=0?}
Clearly $C$ is a cactus and there is a $1$-Lipschitz map $h:C\to X$. We will show that $h$ is in fact a quasi-isometry. We show first that any point in $X$ is at distance less than
$M$ from $h(C)$.

Let $x\in X$ and let $[e,x]$ be a geodesic joining $x$ to $e$. 
Let $x_1$ be the last point (starting from $e$) on this geodesic that lies
in the boundary $\bd Y$ of a graph of (some) level $k$.
(Such $x_1$ exists since $X$ is a graph.) If there is no such point clearly $d(e,x)\leq M$ so $d(x,h(C))\leq M$.

%a connected component of $X\setminus N_M(h_k(C_k))$ for some $k$. 
Otherwise by lemma \ref{close} there is some point $z$ in $h_{k}(C_{k})$ at distance at most $M/10$ from $x_1$.
It follows that
if $d(x,z)\geq M$, since $N_M(h_{k}(C_{k}))$ contains $x_1$, there is some $x_2$ after $x_1$
that lies
in the boundary of a graph of level $k+1$ contradicting our definition of $x_1$.
Therefore $$d(x, h_{k+1}(C_{k+1}))<M.$$

In order to show that $h$ is a quasi-isometry it remains to show that if $x,y\in C$ then $$d(h(x),h(y))\geq cd_C(x,y)-c$$ for some $c>0$.
To do this we introduce first some terminology. We constructed $C$ inductively in stages. In stage $k+1$, to define $C_{k+1}$ we added some circles and points
and we joined each one of them to a circle or point of $C_k$. We called these circles or points \textit{nodes} of the cactus. We will say that the nodes
added at stage $k+1$ are the nodes of \textit{level} $k+1$. We joined each node $S$ of level $k+1$ to a node $S'$ of level $k$ by a segment $p$. We say that $S\cap p$ is the
\textit{basepoint} of the node $S$ and that $S'\cap p$ is an \textit{exit point} of the node $S'$. Of course if $S$ is a point the basepoint of $S$ is equal to $S$. 
We note that $$\length (p)\leq M,$$
so each node at level $k+1$ is at distance at most $M$ from a node of level $k$. It follows that it is enough
to show the inequality $d(h(x),h(y))\geq cd_C(x,y)-c$ for $x,y$ lying in some nodes of $C$, where we denote by $d_C$ the distance in $C$. We remark that each node of level $k>0$ is connected to a unique 
node of level $k-1$.

We have the following easy lemma:

\begin{Lem} \label{fork} If $x,y\in C$ then every geodesic segment joining $x,y$ goes through the same finite set of nodes. 
Any path $\gamma$ from $x$ to $y$ in $C$ goes through this set of nodes.

Moreover there is a node $S$ and points $a,b\in S$ such that any such geodesic contains $\overline {ab}$ (or $S\backslash \overline{ab}$ if $a,b$ are antipodal on $S$), $S$ is the unique minimum level node of the path
and  $\overline {ab}$ is a maximal subarc of $S$ with this property.

\end{Lem}
If the minimal level node $S$ is a point, then $a=b=S$.
\proof Left to the reader.
\qed

If $x,y\in C$ and $S,a,b$ are as in lemma \ref{fork} any geodesic $\gamma $ joining $x,y$ in $C$ can be written
as $$\gamma =\gamma _1\cup \gamma _2 \cup \gamma _3$$ where
$\gamma _1$ is a geodesic joining $x,a$, $\gamma _2$ is a geodesic in $S$ joining $a,b$ and $\gamma _3$ is a geodesic joining $y,b$.
We will show that any path $p$ joining $h(x),h(y)$ in $X$ can be similarly broken in 3 paths, $p_1,p_2,p_3$ have lengths
`comparable' to (or ``longer'' than) those of $\gamma _1, \gamma _2, \gamma _3$.

For $x\in C$, $h(x)$ is a point in $X$, but we may write it as $x$ for
simplicity in the following. 

\begin{Lem} \label{dist-bd} If $Y_2\subseteq Y_1$ are graphs of levels $k+1,k$ respectively then $$d(\bd Y_1,\bd Y_2)\geq 8M/10.$$

\end{Lem}
\proof 
Let $S$ be the node associated to $Y_1$. Then, by the way $Y_2$ was defined, each point of $S$ is at distance $\geq M$ from $\bd Y_2$.
By lemma \ref{close} each point of $\bd Y_1$ 
%\kf{fixed a typo} 
is at distance $\leq M/10+100m<2M/10$ from $S$ so $$d(\bd Y_1,\bd Y_2)\geq 8M/10.$$
\qed

\begin{Lem}\label{2-levels2}
Let $S$ be a node of level $k$ that corresponds to a graph of lebel $k$, $Y_k$.
Then $\diam (S \backslash Y_k) \le 3M$.

\end{Lem}
\proof
If $S$ is a point, then by definition, it is contained
in $\partial Y_k$, so we assume $S$ is a circle. 
Let $S'$ be the node that corresponds to the graph $Y_{k-1}$ with $Y_k \subset Y_{k-1}$
in the construction of $C$. 
We  know that $\diam S \cap Y_k \ge M/20$ by Lemma \ref{contained},
\ref{contained2-2}, depending on $S'$ is a point or a circle. 

First, suppose $S'$ is a circle. Then Lemma \ref{bigpart} applies, and 
as we said in the proof of Lemma \ref{contained2-2}, there
are points $a_2,b_2 \in S$ such that 
$$\overline{a_2qb_2} \subset S\cap Y_k.$$
But since $d(a_2,S')=M$, we have 
$d(a_2,s) \le M + 2000m$
by Lemma \ref{2circles}. Also, $d(b_2,s) \le M+2000m$.
This implies the length of the arc $\overline{a_2sb_2} \le 2M+4000m
\le 3M$, which implies $\diam(S \backslash Y_k) \le 3M$.

Second, suppose $S'$ is a point. The argument is similar
and easier (refer to Remark \ref{contained.remark} instead of Lemma \ref{bigpart}, and use Lemma \ref{geodcirc}
instead of Lemma \ref{2circles}), and we omit it. 
\qed

\begin{Lem}\label{cactus-bd.2}
Let $S_k,S_{k+1}$ be nodes of $C$ corresponding to the graphs $Y_k, Y_{k+1}$ of levels $k,k+1$ respectively
where $Y_{k+1}\subseteq Y_k$.
Let $x\in \partial Y_k$ and  $y\in \partial Y_{k+1}$.
Let $x' \in S_k$ with $d(x,x') \le M/10+100m$ and $y' \in S_{k+1}$
with $d(y,y') \le M/10+100m$.
Then
$$2d(x,y) \ge d_C(x',y') \ge d(x,y)/2$$

\end{Lem}
\proof
There are four cases depending on the two nodes are points or circles.
We only argue that case that both are circles. The other cases 
are similar. 
First, we have $|d(x,y) -d(x',y')| \le 2M/10+200m$.
By Lemma \ref{dist-bd}, $d(x,y) \ge 8M/10$.
Also, by Lemma \ref{2circles},
$|d_C(x',y') -d(x',y')| \le 2000m \le M$.
From those, the conclusion easily follows. 
%\kf{say more?}
\qed

\begin{Lem} \label{basepoint2} Let $S_k,S_{k+1}$ be nodes of $C$ corresponding to the graphs $Y_k, Y_{k+1}$ of levels $k,k+1$ respectively
where $Y_{k+1}\subseteq Y_k$. 
Suppose there is a node $S_{k+1}'$
corresponding to a graph $Z_{k+1}$ of level $k+1$
with $Y_{k+1}\not= Z_{k+1} \subseteq Y_k$.
 Suppose $y\in \bd Y_{k+1}, z\in \bd Z_{k+1} $ 
 and $x \in S_{k+1}, x' \in S_{k+1}'$ with $d(x,y) \le M/10 +100m, 
 d(x',z) \le M/10+100m$ are given. 
 Then
$$d_C(x,x')\leq d(y,z)+5M.$$

\end{Lem}
\proof
Each node $S_k, S_{k+1}, S_{k+1}'$ is a point or a circle, but 
we only discuss the case that all of them are circles. The other cases
are similar (and easier).
Let $a,b\in S_k$ be the point where the nodes $S_{k+1}, S_{k+1}'$
are connected to $S_k$ in $C$.
Then, since $d(y, S_k)=M$, we have $d(x,S_k) \le M+M/10+100m$.
This implies $d_C(a,x) \le M+M/10 + 1100m \le M+2M/10$ by Lemma \ref{2circles}.
Similarly, we have
$d_C(b,x') \le M+ 2M/10$.
Also, since $d(y,S_k)=M$, we have $d(a,y) \le d(a,x)+M/10 \le d_C(a,x)+M/10 \le M+3M/10$. Similarly, $d(b,z) \le M+3M/10$. 
So, 
\begin {align*}
d_C(x, x') & \le d_C(x,a) + d_C(a,b) + d_C(b,x')
\le d_C(a,b) + 2(M+2M/10) 
\\
&=d(a,b)+2(M+2M/10) \le d(a,y)+d(y, z)+d(z,b)
 + 2(M+2M/10) 
 \\
 &\le d(y,z)+2(M+2M/10) + 2(M+3M/10) 
\le d(y, z) +5M
    \end {align*}
\qed

We will show now that if $x,y$ lie in some nodes of $C$ then
$$d(x,y)\geq \dfrac {d_C(x,y)}{2}-20M \  \ \ \ (*).$$
%\kf{swapped  $d$ and $d_C$}
Let's say that $x,y$ lie respectively on nodes $S_k,S_n$ of levels
$k,n$ and the shortest path in $C$ contains an arc $\overline{ab}$ of a (unique) node $S_r$ of minimal level $r$ (and $\overline{ab}$ is maximal with this property) 
by Lemma \ref{fork}. Then we  have that there is a sequence of graphs 
$Y_r,Y_{r+1},...,Y_k$ of levels $r,r+1,...,k$ with $S_r$ corresponding to $Y_r$ and $x\in S_k$. Similarly there is a sequence of graphs $Y_r,Z_{r+1},...,Z_n$  for $y$.

We outline now the proof of the inequality before writing out precise inequalities: For simplicity, we assume that $x \in Y_k$ and 
$y \in Z_n$. Then any path $p$ joining $x,y$ in $X$ intersects
$\bd Y_i$ for $i=r+1,...,k$ by Lemma \ref{fork}.
However by lemma \ref{cactus-bd.2} the part of $p$ from $x$
to $Y_r$ has length comparable to (or longer than)
 the geodesic joining $x$ to $S_r$ in $C$. A similar argument
applies to the corresponding graphs, say $Y_r,Z_{r+1},...,Z_n$  for $y$. Also the intersection points of $p$ with $\bd Y_{r+1},\bd Z_{r+1}$ are `close' to $a,b$ so the
subpath of $p$ joining them has distance comparable to $d_C(a,b)$
by Lemma \ref{basepoint2}.

We give now a detailed argument.

It is possible that $x\in Y_{k-1}$ or $x\in Y_{k-2}$. However in this case by lemma \ref{2-levels2} 
there is some $x_1\in Y_k \cap S_k$ with $d(x,x_1)\leq 3M$.
So we may replace $x$ by $x_1$ and prove $(*)$ with a $20M$ replaced by $15M$, so we assume now that $x\in Y_k$. Similarly we assume that $y\in Z_n$
at a cost of changing the constant of $(*)$ again, say to $10M$.

Let $p$ be a geodesic path joining $x,y$ in $X$. 
Let $$p=p_1\cup p_2\cup p_3$$
where $p_1$ joins $x$ to $\bd Y_{r+1}$, and $p_3$ joins $y$ to $\bd Z_{r+1}$.

\begin{figure}[htbp]
\hspace*{-3.3cm}     
\begin{center}
                                                      
\includegraphics[scale=0.500]{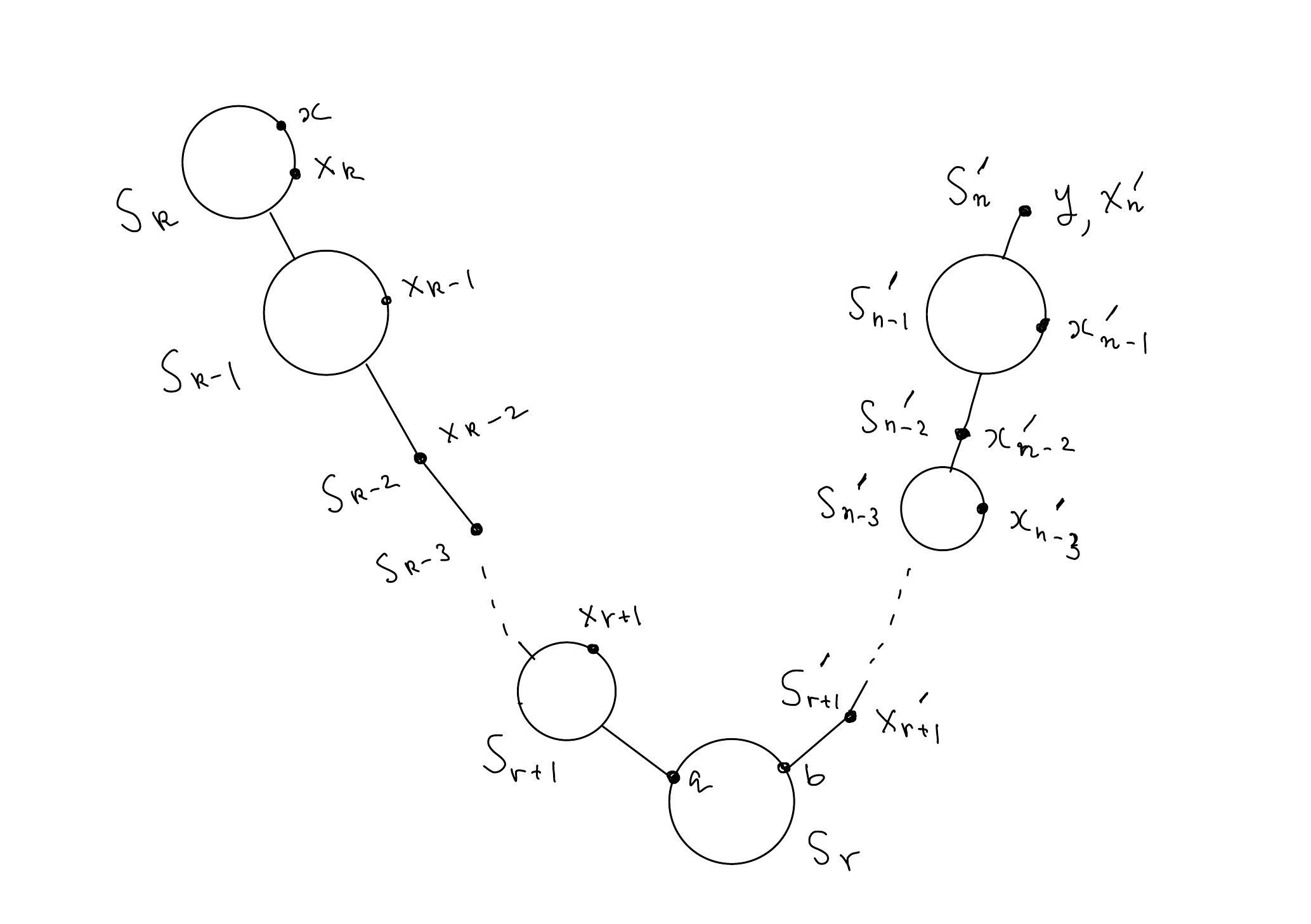}
  \end{center}
\caption{Those points are on $C$, but also in $X$ by the map $h$. Each $x_i$ is almost on $\bd Y_i$ and $x_i'$
is almost on $\bd Z_i$}
\label{estimate.fig}
\end{figure}

The path $p_1$ intersects 
$\bd Y_i$ for $i=r+1,...,k$. 
Let $y_i \in \bd Y_i (r+1 \le i \le k)$ be the first point
that $p_1$ intersects $\bd Y_i$ when it traverse from $\partial Y_{r+1}$ to $x$.
For each $y_i$, take a point $x_i \in S_i$, where $S_i$ is the 
node corresponding to $Y_i$ with $d(y_i,x_i) \le M/10+100m$ (by Lemma \ref{circle} and  \ref{close}).
We obtain the sequence of points $x, x_k, x_{k-1}, \cdots, x_{r+1}$,
which we can see as points on $C$. Join the consecutive two points
by a geodesic in $C$, and we obtain a path, $\beta_1$, from $x$ to $x_{r+1}$ in $C$. See Figure \ref{estimate.fig}.

First, we have
$$d_C(x,x_k) - M/10-100m  =d(x,x_k)-M/10-100m \le d(x,y_k).$$
By Lemma \ref{cactus-bd.2}, 
for each $i, (r+2 \le i \le k)$, we have
$$d_C(x_i, x_{i-1}) \le 2d(y_i,y_{i-1})$$
Summing them up, we have 
\begin{align*}
&d_C(x,x_{r+1}) - M/10 -100m
 \le |\beta_1| -M/10 -100m
 \\
 &\le d_C(x,x_k) -M/10 -100m + d_C(x_k, x_{k-1}) + \cdots + d_C(x_{r+2}, x_{r+1})
 \\
 & \le 2(d(x,y_k) + d(y_k, y_{k-1}) + \cdots + d(y_{r+2},y_{r+1}))=2|p_1|.
 \end{align*}

Similarly, the path $p_3$ intersects
$\bd Z_i$ for $r+1 \le i \le n$.
Let $z_i \in \bd Z_i$ be the first point that $p_3$ intersects $\bd Z_i$
when it traverses to $y$.
For each $z_i$, take a point $x_i' \in S_i'$, where $S_i'$
is the note corresponding to $Z_i$ with $d(z_i,x_i') \le M/10+100m$.
We obtain a sequence of point $y, x_n', x_{n-1}', \cdots, x_{r+1}'$
on $C$, and a path $\beta_3$ in $C$  from $y$ to $x_{r+1}'$
joining the consecutive two points on the sequence by a geodesic
in $C$.
Then as before we have 
$$d_C(y, x_{r+1}') -M/10-100m  \le |\beta_3|-M/10-100m \le 2|p_3|. $$

Finally, by Lemma \ref{basepoint2}
$$ d_C(x_{r+1}, x_{r+1}') - 5M \le d(y_{r+1}, z_{r+1}) =|p_2|. $$
Combining them, we have

\begin{align*}
2d(x,y)&=2(|p_1|+|p_2|+|p_3|) 
\\
&\ge d_C(x,x_{r+1})+d_C(x_{r+1},x_{r+1}')
+ d_C(x_{r+1}',y) -5M-2M/10-200m 
\\
&\ge d_C(x,y)-6M.
\end{align*}
We got $d(x,y) \ge d_C(x,y)/2 - 3M$. 
The inequality (*) is shown. 
The quasi-isometry constants depend only on $M$, so that only
ono $m$.
The proof is complete. 
\qed

\begin{Cor}\label{bottlecor}
Let $X$ be a geodesic metric space.
$X$ is quasi-isometric to a cactus
if and only if it satisfies $\sharp $.

Moreover the quasi-isometry-constants depend only on $m$ in $\sharp$.
%\kf{added on uniformness}
\end{Cor}
\proof Every metric space is quasi-isometric, say, $(2,1)$-quasi-isometric,
 to a graph and $\sharp $ is clearly a condition that is invariant 
under quasi-isometries.
Uniformity is clear. 
\qed

\section{Characterization of cacti by fat theta curves}\label{sec-char-proof}

We will prove Lemma \ref{equivalence}, Lemma \ref{fatissharp} and 
Proposition \ref{prop-char-cactus}. As usual we only argue for the case that $X$ is a graph. 
First we prove the two lemmas. 
\\
\\
{\it Proof of Lemma \ref{equivalence}}.
If $X$ satisfies $(\sharp, m)$ for $m>0$, then 
$X$ is quasi-isometric to a cactus, $C$, and the quasi-isometry constants
depend only on $m$ by Corollary \ref{main}.
We may assume the quasi-isometry is continuous, and suppose 
it is a $(K,L)$-quasi-isometry, where $K,L$ depend only on $m$. 

Set $M=K+KL$, which depends only on $m$. Then $X$ contains no $M$-fat theta curve, since 
if there was 
an $M$-fat theta curve in $X$, then the quasi-isometry would give 
a $1$-fat theta curve in $C$. But the cactus does not contain 
any embedded fat theta curve, a contradiction.  
\qed

It would be interesting to find an elementary proof of Lemma  \ref{equivalence}
without using 
Corollary \ref{main}.\
%kf {added a comment}
\\
\\
{\it Proof of Lemma \ref{fatissharp}}.
In the proof below we don't try to optimize the relationship between $m,M$, we simply pick our constants so that all inequalities we need are easy to verify.

We prove the contrapositive. Let $m=1000M$.  Assume that there are two points $a,b$ at distance $\geq 10m$ such that if  $x,y$
are at distance $\geq 4m$ from both $a,b$ then there is a path in $$X\setminus (B_m(x)\cup B_m(y))$$ joining $a,b$.

Let $\gamma $ be a geodesic joining $a,b$. We say that a path $B$ is an $m/10$-bridge if
$B$ is a union of paths $L,D,R$ such that $L=[l_1,l_2], R=[r_1,r_2]$ are geodesic paths of length $[m/10]$
with $l_2,r_2\in \gamma $ and $d(l_2,r_2)\geq m$, $D$ has endpoints $l_1,r_1$ and every $x\in D$ satisfies
$$d(x,\gamma )\geq [m/10].$$
Here, $[x]$ means the floor function.
%\kf{made a comment on $[x]$}

By our assumption it follows $m/10$ bridges exist.
To see it, let $c$ be the mid point of $\gamma$. Denote $\gamma_1 \subset \gamma$
be the part between $a,c,$ and set $\gamma_2=\gamma \backslash \gamma_1$. 
Then there is a path, $\alpha$, between $a,b$ in $X \backslash B_m(c)$
by our assumption. 
Let $\alpha(t_1)$ be the last point on $\alpha$ that is in $N_{[m/10]}(\gamma_1)$, and 
$\alpha(t_2)$ be the first point that is in $N_{[m/10]}(\gamma_2)$
with $t_1 < t_2$. Such $t_1, t_2$ must exist. 
Let $\alpha_0$ be the part of $\alpha$ between $\alpha(t_1), \alpha(t_2)$.
Let $c_1\in \gamma$ be a point with $d(\alpha(t_1), c_1)=[m/10]$, and 
$c_2 \in \gamma$ a point with $d(\alpha(t_2),c_2)=[m/10]$.
Then 
$$[\alpha(t_1),c_1) \cup \alpha_0 \cup [\alpha(t_2),c_2)]$$
is an $m/10$ bridge. 
%\kf{added a construction of a bridge}

 Let $$B=[l_1,l_2]\cup D\cup [r_1,r_2]$$ be an $m/10$ bridge for 
which $d(l_2,r_2)$ is maximal. We modify now the bridge $B$ as follows. If $x=D(t)$ is a point of $D$ that is at distance $\geq m/3$
from $l_1,r_1$ and $y=D(s)$ is a point such that $d(D(s),D(t))\leq m/100$ and $|s-t|>m/100$
we replace $D([s,t])$ by a geodesic path with the same end-points. Here we assume as usual that $D$ is parametrized by arc length.

We still call $D$ the new path
we obtain after this replacement.  We keep doing this operation as long as
$d(D(t),\gamma )\geq m/100$ for all $t$. It is clear that one can do this operation finitely many times.

There are two cases.\smallskip

\textit{Case 1}. After we do this operation finitely many times
there is a $t$ such that $d(D(t),\gamma )< m/50$
for some $t$. We set then $z_1=D(t)$ and we pick a closest point $z_2\in \gamma $ ie a point such that
$$d(z_1,z_2)=d(z_1,\gamma ).$$
We claim that we have now an $m/1000$-fat theta curve obtained as follows: The vertices of the theta are $z_1,z_2$. The three arcs
of the theta that are $m/1000$ away are $[l_1,l_2], [z_1,z_2], [r_1,r_2]$ and the endpoints of these arcs are joined by the obvious
subarcs of $D$ and $\gamma $, which gives two subsets that are $m/500$-away.  \smallskip

\textit{Case 2}. After we do this operation finitely many times we finally obtain an arc $D$ such that $d(D(t),\gamma )\geq m/100$
for all $t$ and if $x=D(t)$ is a point of $D$ that is at distance $\geq m/3$
from $l_1,r_1$ and $d(D(s),D(t))\leq m/100$ then $D([s,t])$ is a geodesic path. We pick now $z_1$ on $D$ such that is at distance $\geq m/3$
from $l_1,r_1$ and $z_2$ a point on $\gamma $ between $l_2,r_2$ and at distance $\geq m/3$
from $l_2,r_2$. By our assumption there is a path $\beta $ in
$$X\setminus (B_m(z_1)\cup B_m(z_2))$$ joining $a,b$. 

If we write $\gamma $ as union of successive geodesics $\gamma =\gamma_1\cup [l_2,r_2]\cup \gamma _3$ we will define a theta curve with vertices $a,b$. Two of the arcs of the theta curve are
the curves: $\gamma, \eta=\gamma _1\cup [l_2,l_1]\cup D\cup [r_1,r_2]\cup \gamma _2,$.

Let $p _1,p _2$ to be geodesic subpaths of length $m/100$ centered
at $z_1,z_2$ respectively (so they are subarcs of $\gamma, \eta$. Let's say that $z_1=\eta (s_1),z_2=\gamma (s_2)$.

Let $$u_1=\max \{t:d(\beta (t),\eta ([0,s_1])\cup \gamma ([0,s_2]))\leq m/1000\}$$
If $\length (\eta)=\ell _1, \length (\gamma)=\ell _2$ we set
$$u_2=\min \{t: t>u_1 \text{ and } d(\beta (t),\eta ([s_1,\ell_1])\cup \gamma ([s_2,\ell_2]))\leq m/1000\}.$$

We pick shortest paths: $\beta _1=[\beta (u_1),b_1]$ joining $\beta (u_1)$ to $\eta ([0,s_1])\cup \gamma ([0,s_2])$
and $\beta _2 =[\beta (u_2),b_2]$ joining $\beta (u_2)$ to $\eta ([s_1,\ell_1])\cup \gamma ([s_2,\ell_2])$. We finally pick arcs: $\beta _3\subseteq \eta ([0,s_1])\cup \gamma ([0,s_2])$
joining $a$ to $b_1$ and $\beta _4\subseteq \eta ([s_1,\ell_1])\cup \gamma ([s_2,\ell_2])$
joining $b_2$ to $b$. Now the third arc of our theta curve is the arc
$$\beta '=\beta _3\cup \beta _1\cup \beta ([u_1,u_2])\cup \beta _2\cup \beta _4$$
and we define the subarc $p_3$ to be
$$p_3=\beta ([u_1,u_2]).$$

We note that the distance between $\eta ([0,s_1-m/100])\cup \gamma ([0,s_2-m/100])$ and $\eta ([s_1+m/100,\ell_1])\cup \gamma ([s_2+m/100,\ell_2])$ at least $m/100$.

It follows from this and our definition of $p_3$ that the arcs $p_1,p_2,p_3$ satisfy the conditions of the definition
of an $M=m/1000$ fat theta curve.
\qed
\\
The proposition immediately follows. 
\\
\\
{\it Proof of Proposition \ref{prop-char-cactus}}.
We combine Lemma \ref{equivalence} and Lemma \ref{fatissharp}.
\qed

We conclude with a lemma that confirms the speculation from \cite{FP} we mentioned in the 
introduction.
\begin{Lem}\label{lem-annuli}
Let $P$ be a plane with a geodesic metric with a
base point $e$, and let $A=A(r,r+m)$ be the 
set of points $x$ with $r \le d(e,x) < r+m$.

Then each connected component of  $A=A(r,r+m)$
with the path metric satisfies $(\sharp, m)$.
Therefore, it is  uniformly quasi-isometric to a cactus with 
the quasi-isometric constants depending only on $m$. 

\end{Lem}
\proof
In \cite[Lemma 4.1]{FP} we proved that a connected component of $(A,d_A)$
does not contain any $m$-fat theta curves. 
Then by Lemma \ref{fatissharp} it satisfies $(\sharp, M)$ for 
$M=1000m$ (see the proof of the lemma too).
Then by Corollary \ref{main} it is uniformly quasi-isometric to a cactus
and the quasi-isomety constants depend only on $m$. 
\qed

\end{document}